\documentclass[11pt,leqno,a4paper]{article}

\usepackage[hypertex]{hyperref}
\usepackage{amsfonts}
\usepackage{latexsym}
\usepackage{amssymb,amsmath}
\usepackage{epsfig}

\newcommand{\C}{\mathbb C}
\newcommand{\R}{\mathbb R}
\newcommand{\N}{\mathbb N}

\newcommand{\be}{\begin{equation}}
\newcommand{\ee}{\end{equation}}
\newcommand{\ba}{\begin{eqnarray}}
\newcommand{\ea}{\end{eqnarray}}
\newcommand{\eq}[1]{\eqref{#1}}

\def\qd{{\rule{2.3mm}{2.3mm}}}
\def\qed{{\hfill$\quad$\qd\\}}

\def\ep{{\varepsilon}}

\newtheorem{thm}{\bf Theorem}
\newtheorem{lem}[thm]{\bf Lemma}

\newtheorem{rem}[thm]{\bf Remark}

\renewcommand{\leq}{\leqslant}
\renewcommand{\geq}{\geqslant}
\renewcommand{\le}{\leqslant}
\renewcommand{\ge}{\geqslant}

\date{January 28, 2010}

\begin{document}

\title{Null controllability of a parabolic system with a cubic coupling term}

\author{
Jean-Michel Coron \thanks{Institut universitaire de France and Universit\'{e} Pierre et Marie Curie - Paris 6,
UMR 7598 Laboratoire Jacques-Louis Lions, Paris, F-75005 France
({\tt coron@ann.jussieu.fr}).  JMC was partially
supported by the ``Agence Nationale de la Recherche'' (ANR), Project C-QUID,
grant BLAN-3-139579.}
\and Sergio Guerrero \thanks{Universit\'{e} Pierre et Marie Curie - Paris 6,
UMR 7598 Laboratoire Jacques-Louis Lions, Paris, F-75005 France
({\tt guerrero@ann.jussieu.fr}). SG was partially
supported by the ``Agence Nationale de la Recherche'' (ANR), Project CISIFS,
grant ANR-09-BLAN-0213-02.}
\and Lionel Rosier \thanks{Institut \'{E}lie Cartan,
UMR 7502 UHP/CNRS/INRIA,
B.P. 239,  54506 Vand\oe uvre-l\`{e}s-Nancy Cedex, France.
({\tt rosier@iecn.u-nancy.fr}). LR was partially
supported by the ``Agence Nationale de la Recherche'' (ANR), Project CISIFS,
grant ANR-09-BLAN-0213-02.}}
\maketitle



\begin{abstract}
We consider a system of two parabolic equations
with a forcing control term present in one equation
and a cubic coupling term in the other one. We prove that the system
is locally null controllable. 

\end{abstract}
\section{Introduction}
\setcounter{equation}{0}
The control of coupled parabolic systems is a challenging
issue, which has attracted the interest of the control community
in the last decade. Let $\Omega$ be a nonempty connected bounded subset of $\R^N$ of class $C^2$. Let
$\omega$ be a nonempty open subset of $\Omega$. In \cite{ABDG2007} and \cite{ABDG2009}, the authors identified sharp conditions for the control of systems of the form
\begin{eqnarray}
w_t&=& (D\Delta +A)w + B  h,
\end{eqnarray}
where $w=(w_1,...,w_N):\Omega \rightarrow \R^N$ is the state to be controlled,
$h=h(t,\cdot):\Omega \rightarrow \R^M$ is the control input supported in
$\omega$, and $D: \R ^N \rightarrow \R ^N$ is a diagonal operator, $A:\R ^N\rightarrow \R ^N$, and
$B:  \R^M \rightarrow \R^N$ are linear maps. In general, the rank of $B$ is less
than $N$, so that the controllability of the full system depends strongly on the
(linear) coupling present in the system.
See \cite{2000-Teresa,Gonzalez-Burgos-Perez-Garcia, SG1} for related results. See also \cite{2009-FGT} for boundary controls, \cite{Assiaetal} for some inverse problems and \cite{2006-Fernandez-Cara-Guerrero-etal,SG2,2009-Coron-Guerrerero} for the Stokes system.

Here, we are concerned with the control of semilinear parabolic systems in which the coupling occurs through {\em nonlinear} terms only. More precisely, we
study the control properties of systems of the form
\begin{equation}\label{system-0}
    \left\{\begin{array}{ll}
    \displaystyle u_{t}-\Delta u=g(u,v)+h1_{\omega}&\hbox{ in }(0,T)\times \Omega, \\ \noalign{\medskip} v_{t}-\Delta
    v=u^3+ Rv
        &\hbox{ in } (0,T)\times\Omega,\\ \noalign{\medskip}
         u=0,\quad v=0&\hbox{ on }(0,T)\times\partial\Omega,
    \end{array}\right.
\end{equation}
where $g:\mathbb{R}\times\mathbb{R}\rightarrow\mathbb{R}$
is a given function of class $C^\infty$ vanishing at $(0,0)\in \mathbb{R}\times\mathbb{R}$, $R$ is a given real number and $1_{\omega}$ is the characteristic function of $\omega$. This a control system where, at time $t\in [0,T]$, the state is $(u(t,\cdot), v(t,\cdot)): \Omega \rightarrow \R^2$ and the control is $h(t,\cdot): \Omega\rightarrow \R$.

The goal of this paper is to prove the local null controllability of system (\ref{system-0}). Our main result is as follows.
\begin{thm}\label{th}
There exists $\delta>0$
such that, for every $(u_0,v_0)\in L^\infty(\Omega)^2$ satisfying
$$
\|u_{0}\|_{L^\infty(\Omega)}+\|v_{0}\|_{L^\infty(\Omega)}<\delta,
$$
there exists a control $h\in L^{\infty}((0,T)\times \Omega)$ such that the  solution $(u,v)\in
L^\infty([0,T]\times\Omega)^2$ of the Cauchy problem
\begin{equation}
\label{system}
    \left\{\begin{array}{ll}
    \displaystyle u_{t}-\Delta u=g(u,v)+h1_{\omega}&\hbox{ in }(0,T)\times \Omega, \\ \noalign{\medskip} v_{t}-\Delta
    v=u^3+Rv
        &\hbox{ in } (0,T)\times\Omega,\\ \noalign{\medskip}
         u=0,\quad v=0&\hbox{ on }(0,T)\times\partial\Omega,
         \\ \noalign{\medskip}
u(0,\cdot)=u_{0}(\cdot),\quad v(0,\cdot)=v_{0}(\cdot)&\hbox{
in }\Omega,
    \end{array}\right.
\end{equation} satisfies
\begin{equation}\label{condT}
u(T,\cdot)=0\quad\hbox{and}\quad
v(T,\cdot)=0\hbox{\,\, in }\Omega.
\end{equation}
\end{thm}

Let us give a system from Chemistry to which our result applies.
A reaction-diffusion system describing a reversible
chemical reaction (see \cite{BH2003, CHKO, ET1989})
takes the form
\begin{eqnarray}
u_t &=& \Delta u -a k(u^k-v^m)  \qquad\text{ in } (0,T)\times \Omega \label{CR1},\\
v_t &=& \Delta v + b k(u^k-v^m) \qquad\text{ in } (0,T)\times \Omega \label{CR2},
\end{eqnarray}
together with homogeneous Neumann boundary conditions.
In \eqref{CR1}-\eqref{CR2}, $a$ and $b$ denote some positive numbers, and $k$ and $m$ are positive integers. The corresponding reversible chemical reaction reads $kA \rightleftharpoons mB$. Incorporating a forcing term $1_\omega h$ in \eqref{CR1},
we obtain a system of the form \eqref{system} when $k=3$ and $m=1$, so that
Theorem \ref{th}
may be applied. (In fact, Theorem \ref{th} deals with Dirichlet homogeneous boundary conditions, but the proof we give here can easily be adapted to deal with homogeneous Neumann boundary conditions as well, and a scaling argument shows that one may assume without loss of generality that $b=1/3$).

\begin{rem}
In Theorem \ref{th}, it is not possible to replace $u^3$ by $u^2$. Indeed, by the maximum principle, for every $(u_0,v_0)\in L^\infty(\Omega)^2$ and for every $h\in L^{\infty}((0,T)\times \Omega)$, the  solution $(u,v)\in
L^\infty([0,T]\times\Omega)^2$ of the Cauchy problem
\begin{equation*}
    \left\{\begin{array}{ll}
    \displaystyle u_{t}-\Delta u=g(u,v)+h1_{\omega}&\hbox{ in }(0,T)\times \Omega, \\ \noalign{\medskip} v_{t}-\Delta
    v=u^2+Rv
        &\hbox{ in } (0,T)\times\Omega,\\ \noalign{\medskip}
         u=0,\quad v=0&\hbox{ on }(0,T)\times\partial\Omega,
         \\ \noalign{\medskip}
u(0,\cdot)=u_{0}(\cdot),\quad v(0,\cdot)=v_{0}(\cdot)&\hbox{
in }\Omega,
    \end{array}\right.
\end{equation*}
if it exists, satisfies
$$
v(T,\cdot )\geqslant v^*(T,\cdot) \text{ in } \Omega,
$$
where $v^*\in L^\infty([0,T]\times \Omega)$ is the solution of the linear Cauchy problem
\begin{equation*}
    \left\{\begin{array}{ll}
   v_{t}^*-\Delta
    v^*=Rv^*
        &\hbox{ in } (0,T)\times\Omega,\\ \noalign{\medskip}
         v^*=0&\hbox{ on }(0,T)\times\partial\Omega,
         \\ \noalign{\medskip}
v^*(0,\cdot)=v_{0}(\cdot)&\hbox{
in }\Omega .
    \end{array}\right.
\end{equation*}
In particular, by the (strong) maximum principle, if $v_0\geqslant 0$ and $ v_0\not = 0$, then $v(T,\cdot ) >0$ in $\Omega$.
L. Robbiano asked to the authors whether the result in Theorem \ref{th},
still with $u^3$ replaced by $u^2$, could be true if we consider  {\em complex-valued} functions. The following result,
whose proof is sketched in Appendix, shows that this is indeed the case.
\end{rem}
\begin{thm}\label{th2}
There exists $\delta>0$
such that, for every $(u_0,v_0)\in L^\infty(\Omega ;\C)^2$ satisfying
$$
\|u_{0}\|_{L^\infty(\Omega)}+\|v_{0}\|_{L^\infty(\Omega)}<\delta,
$$
there exists a control $h\in L^{\infty}((0,T)\times \Omega;\C)$ such that the  solution $(u,v)\in
L^\infty([0,T]\times\Omega;\C)^2$ of the Cauchy problem
\begin{equation}
\label{systemth2}
    \left\{\begin{array}{ll}
    \displaystyle u_{t}-\Delta u=g(u,v)+h1_{\omega}&\hbox{ in }(0,T)\times \Omega, \\ \noalign{\medskip} v_{t}-\Delta
    v=u^2+Rv
        &\hbox{ in } (0,T)\times\Omega,\\ \noalign{\medskip}
         u=0,\quad v=0&\hbox{ on }(0,T)\times\partial\Omega,
         \\ \noalign{\medskip}
u(0,\cdot)=u_{0}(\cdot),\quad v(0,\cdot)=v_{0}(\cdot)&\hbox{
in }\Omega,
    \end{array}\right.
\end{equation} satisfies
\begin{equation}\label{condT2}
u(T,\cdot)=0\quad\hbox{and}\quad
v(T,\cdot)=0\hbox{\,\, in }\Omega.
\end{equation}
\end{thm}

When trying to prove a local null controllability result, the first thing to do is to look at the null controllability of the linearized control system around $0$. Here, the linearized control system reads
\begin{equation}\label{system-lin}
    \left\{\begin{array}{ll}
    \displaystyle u_{t}-\Delta u=\partial_u g(0,0)u +\partial_v g(0,0)v+h1_{\omega}&\hbox{ in }(0,T)\times \Omega, \\ \noalign{\medskip} v_{t}-\Delta
    v= Rv
        &\hbox{ in } (0,T)\times\Omega,\\ \noalign{\medskip}
         u=0,\quad v=0&\hbox{ on }(0,T)\times\partial\Omega.

    \end{array}\right.
\end{equation}
Clearly the control $h$ has no influence on $v$ and, if $v(0,\cdot)\not =0$, then $v(T,\cdot)\not =0$. Hence
the linearized control system \eq{system-lin} is not null controllable and this strategy cannot be applied to prove Theorem~\ref{th}.

Our proof of Theorem \ref{th} relies on the return method, a method introduced in \cite{1992-mcss-coron} for a stabilization problem and in \cite{1996-jmpa-coron} for the controllability of the Euler equations of incompressible fluids (see \cite[Chapter 6]{coron-book} and the references therein for other applications of this method). Applied to the control system \eq{system-0}, it consists in
 looking for a trajectory
$(( \overline u,  \overline v),  \overline h)$ of the control system \eq{system-0} such that
\begin{enumerate}
\item it goes from $(0,0)$ to $(0,0)$, i.e. $\overline  u(0,\cdot)=\overline v(0,\cdot)
=\overline u(T,\cdot)=\overline v(T,\cdot)=0$;
\item the linearized control system around that trajectory is null controllable.
\end{enumerate}
With this trajectory and a suitable fixed point theory at hand, one can hope to get
the null controllability stated in Theorem \ref{th}. We shall see that this is indeed the case.

In a forthcoming paper \cite{CGR2011}, we investigate the case of more general nonlinear coupling terms. In particular, this result can be applied to the chemical reaction system (\ref{CR1})-(\ref{CR2}) for any pair $(k,m)$ with $k$ an odd integer  and also to the internal control of the Ginzburg-Landau equation with a control input taking real values. See \cite{Fu2006}, \cite{RZ2009} for the control of the Ginzburg-Landau equation with a complex control input.

The paper is organized as follows. Section \ref{sec-construction-baruetc}
is devoted to the construction of the trajectory $(( \overline u, \overline v), \overline h)$. In section \ref{seclocalnulbarubarv}, using some Carleman inequality, we prove that the linearized control system around $((\overline u, \overline v), \overline h)$ is null controllable (sub-section \ref{subseccontrolabbility-linearized}). Next, we deduce the local null controllability around this trajectory by using the Kakutani fixed-point theorem (sub-section \ref{sec-fixed-point}). The appendix contains a sketch of the proof
of Theorem \ref{th2}.

\section{Construction of the trajectory $((\overline u, \overline v), \overline h)$}
\label{sec-construction-baruetc}

\setcounter{equation}{0}
Let us define $Q:=(0,T)\times \Omega$. The goal of this section is to prove the existence of $\overline u\in C^\infty(Q)$, $\overline v\in C^\infty(Q)$ and $\overline h\in C^\infty(Q)$ such that
\begin{gather}
\label{supportuvh}
\text{the supports of $\overline u$, $\overline v$ and $\overline h $ are compact and  included in $(0,T)\times \omega$,}
\\
\overline u_t-\Delta \overline u= g(\overline u,\overline v)+\overline h
\qquad \text{ in } Q,
\\
\overline v_t-\Delta \overline v= \overline u^3 +R \overline v
\ \ \ \qquad \text{ in } Q,
\\
\label{barunot0}
\overline u \not\equiv 0.
\end{gather}

The existence of such $\overline u\in C^\infty(Q)$, $\overline v\in C^\infty(Q)$ and
$\overline h\in C^\infty(Q)$ follows from the following theorem.
\begin{thm}
\label{thm1}
Let $\rho >0$ and $R$ be two constants.
Then there exist two functions $V: (t,x)\in \R\times \R^N\mapsto V(t,x) \in \R $
and $K: (t,x)\in \R\times \R^N\mapsto K(t,x) \in \R $ such that
\begin{gather}
V\in C^\infty (\R\times \R ^N), \text{ and }
V(t,x)=0\,   \text{ for all } (t,x)\in \R\times\R^N \text{ with } \max(|t|,
|x|)\ge \rho,
\label{E1}
\\
K\in C^\infty (\R\times \R ^N),\text{ and }
K(t,x)=0 \,\text{ for all } (t,x)\in \R\times\R^N \text{ with } \max(|t|,
|x|)\ge \rho,
\label{regk-support}
\\
K\not\equiv 0,
\label{knot9}
\\
V_t=\Delta V + RV + K^3.
\label{E2}
\end{gather}
\end{thm}
Indeed, let $x_0 \in \omega$ and let $\rho>0$ be small enough
 so that
$$
\max(|t-\frac{T}{2}|, |x-x_0|
)\leqslant \rho \Rightarrow (t,x)\in (0,T)\times\omega.
$$
Then, it suffices to define $\overline u\in C^\infty(Q)$, $\overline v\in C^\infty(Q)$ and
$\overline h\in C^\infty(Q)$ by
\begin{gather*}
\overline u(t,x):= K(t-\frac{T}{2},x-x_0), \, \overline v(t,x):= V(t-\frac{T}{2},x-x_0), \,\forall (t,x)\in Q,
\\
\overline h := \overline u_t-\Delta \overline u -g(\overline u , \overline v) \text{ in } Q.
\end{gather*}

\begin{figure}[hbtp]
\begin{center}
\epsfig{file=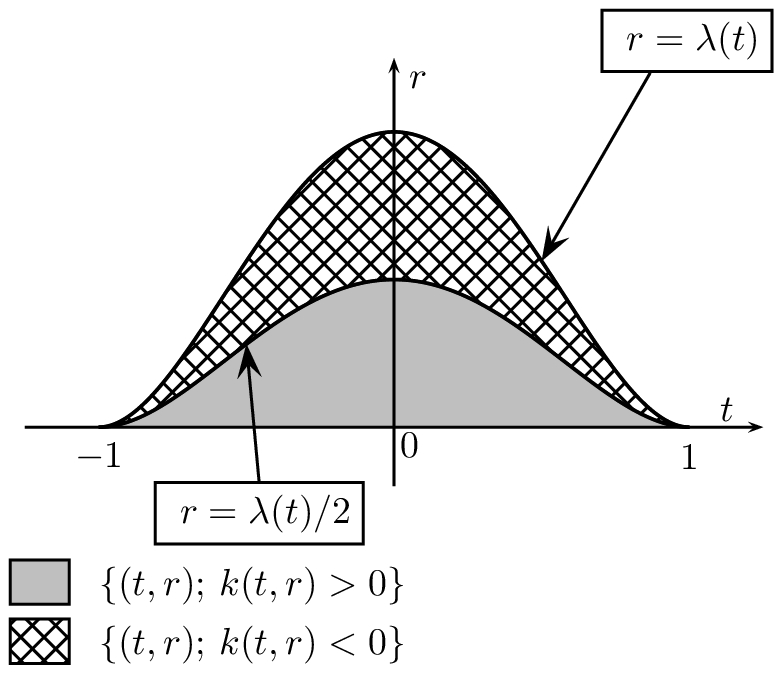, width=10cm}
\caption{$\left\{(t, r);\, k(t,r)>0\right\}$ and
$\left\{(t, r);\, k(t,r)<0\right\}$}
\label{figk>0k<0}
\end{center}
\end{figure}

\textbf{Proof of Theorem \ref{thm1}.}
Note first that we may assume
without loss of generality that $R=0$. Indeed, setting
$$
\tilde V(t,x)=e^{-Rt}V(t,x),\,
\tilde K(t,x)=e^{-Rt/3}K(t,x),
$$
then \eqref{E2} is transformed into
$$
{\tilde V}_t=\Delta\tilde V +\tilde K^3.
$$
\, From now on, we assume that $R=0$. We may also assume that $\rho =1$.
Indeed, if the construction has been done for $\rho =1$
and $R=0$, then
for any $\rho >0$ the functions
$$
\tilde V (t,x)=V(\rho ^{-2}t, \rho ^{-1}x), \quad
\tilde  K (t,x)=\rho ^{-\frac{2}{3}}
K(\rho ^{-2}t, \rho ^{-1}x),
$$
with support in $[-\rho ^2, \rho^2]\times \{ |x| \le \rho \}$,
satisfy the equation
${\tilde V}_t = \Delta \tilde V + {\tilde K}^3$.
We assume from now on that $R=0$ and that $\rho =1$.
Let $r=|x|$. We seek for a radial function $V(t,x)=v(t,r)$
fulfilling the following properties
\begin{eqnarray}
&& v\in C^\infty(\R \times \R ^+),\,  v(t,r)=0 \text{ for }
|t|\ge 1 \text{ or } r\ge 1, \label{E1bis}\\
&& k:=(v_t-v_{rr}-\frac{N-1}{r}v_r)^{\frac{1}{3}}
\in C^\infty (\R \times \R ^+). \label{E2bis}
\end{eqnarray}
The smoothness of $V$ and $K:=(V_t-\Delta V)^{\frac{1}{3}}$
at the points $(t,0)$, $t\in [-1,1]$ will follow from
additional properties of $v$ (see below). As far as
 the construction of $v$ is concerned, the idea is to have a precise
knowledge  of the place where $k$ vanishes, and a good ``behavior'' of $v$ near the place where $k$ vanishes to ensure that $k$ is of class $C^{\infty}$.
For the function $v$ we are going to construct, we shall have
\begin{gather*}
\left\{(t, r);\, k(t,r) < 0\right\}=\left\{(t,r);\ 0<\lambda(t)/2<r<\lambda (t)\right\},
\\
\left\{(t, r);\, k(t,r) > 0\right\}=\left\{(t,r);\ 0<r<\lambda(t)/2\right\}.
\end{gather*}
See Figure \ref{figk>0k<0}.

Let us introduce a few notations. Let
\begin{eqnarray}
\lambda (t) &=&
\ep (1-t^2)^2 , \,  |t|<1, \label{B5}\\
f_0(t)&=&
\left\{
\begin{array}{ll}
e^{-\frac{1}{1-t^2}}, &|t|<1, \\
0,\,  &|t|\ge 1,
\end{array}
\right.
\label{B6}
\end{eqnarray}
where $\ep >0$ is a (small) parameter chosen later.
We search $v$ in the form
\be
v(t,r)=\sum_{i=0}^3 f_i(t)g_i(z),
\label{B1}
\ee
where $z:=r/\lambda (t)$, $g_0$ is defined in Lemma \ref{lem2.1}
(see below), and the functions
$f_i=f_i(t)$, $1\le i\le 3$, and $g_i=g_i(z)$, $1\le i \le 3$,
defined during the proof, are in $C^\infty (\R )$ and fulfill
\begin{eqnarray}
&&\text{supp } f_i\subset [-1,1], \label{B3} \\
&&\text{supp } g_i\subset [\frac{1}{2} -\frac{\delta }{2},
\frac{1}{2} +\frac{\delta}{2}]. \label{B4}
\end{eqnarray}
In \eqref{B4}, $\delta \in (0,1/10)$ is a given number.
Let us begin with the construction of $g_0$.
\begin{lem}
\label{lem2.1}
There exists a function $G\in C^\infty ([0,+\infty ))$
such that
\begin{eqnarray}
G(z) = (z-\frac{1}{2})^3 &\text{ for }&
\frac{1}{2} -\delta <z<\frac{1}{2} +\delta, \label{A1}\\
(z-\frac{1}{2})G(z)>0  &\text{ for }& 0<z<1,\ z\ne \frac{1}{2},
 \label{A2}
\end{eqnarray}
and such that the solution $g_0$ to the Cauchy problem
\begin{eqnarray}
&&{g_0}''(z) +\frac{N-1}{z}{g_0}'(z)= G(z),\quad z>0,
\label{A4}\\
&&g_0(1)=g_0'(1)=0, \label{A5}
\end{eqnarray}
satisfies
\begin{eqnarray}
g_0(z)&=& 1-z^2\text{ if }\  0<z<\delta, \label{A6a}\\
g_0(z)&=& e^{-\frac{1}{1-z^2}}
\text{ if }\ 1-\delta < z <1, \label{A6b}\\
g_0(z)&=& 0 \text{ if }\ z\ge 1. \label{A6c}
\end{eqnarray}
\end{lem}
\textbf{Proof of Lemma \ref{lem2.1}.} Note first that by \eqref{A4}, \eqref{A6a} to \eqref{A6c}, we have
\be
G(z)=
\left\{
\begin{array}{ll}
-2N \quad &\text{ if } 0 < z < \delta , \\
 \left[-2N(1-z^2)^{-2}-8z^2(1-z^2)^{-3}+4z^2(1-z^2)^{-4}\right]
e^{-\frac{1}{1-z^2}}
\quad &\text{ if } 1-\delta <z<1,\\
0 &\text{ if } z\ge 1.
\end{array}
\right.
\label{A7}
\ee
and hence only the values of $G$ on $[\delta, \frac{1}{2}-\delta]$ and
on $[\frac{1}{2} + \delta , 1-\delta ]$ remain to be defined.
Let $ G\in C^\infty (0,+\infty )$ be any function satisfying
\eqref{A1} and \eqref{A7}, and denote by $g_0$ the solution
of \eqref{A4}-\eqref{A5}. Clearly, \eqref{A6b}-\eqref{A6c}
are satisfied. Finally, it is clear that \eqref{A6a} holds
if and only if $g_0(0^+)=1$ and $g_0'(0^+)=0$. Note that
\eqref{A4} may be written as follows:
\be
\label{A10}
\frac{1}{z^{N-1}} (z^{N-1}g_0')'=G.
\ee
Using \eqref{A5}, this gives upon integration
\be
\label{A11}
-z^{N-1} g_0'(z)=\int_z^1 s^{N-1}G(s)\,ds.
\ee
This imposes the condition
\be
\label{A12}
\int_0^1 s^{N-1}G(s)\,ds =0.
\ee
Note that, if \eqref{A12} holds, then, by \eqref{A11},
\be
\label{A13}
g_0'(z)=\frac{1}{z^{N-1}}\int_0^z s^{N-1}G(s)\,ds,
\ee
which, combined to \eqref{A7}, yields
$$
g_0'(z)=-2z,\, 0<z<\delta,
$$
and $g_0'(0^+)=0$. Integrating \eqref{A13} on $[z,1]$ and using (\ref{A5}),
(\ref{A7}), (\ref{A12}) and an integration by parts, we obtain,
 for $0<z<\delta$,
\begin{eqnarray}
g_0(z)
&=&-\int_z^1 \frac{1}{y^{N-1}}(\int_0^y s^{N-1}G(s)\,ds)dy\\
&=&
\left\{
\begin{array}{ll}
\displaystyle\frac{1}{2-N}\int_z^1 yG(y)dy
-\frac{2}{2-N}z^2\quad
&\text{ if } N\ne 2, \\[5mm]
\displaystyle\int_z^1 y(\ln y) G(y) dy -2z^2\ln z\quad
&\text{ if } N=2.
\end{array}
\right.
\label{A14}
\end{eqnarray}
Then $g_0(0^+)=1$ provided that
\begin{equation}
\left\{
\begin{array}{ll}
\displaystyle \int_0^1 y G(y)dy =
2-N   &\text{ if } N\ne 2,\\[5mm]
\displaystyle \int_0^1 y (\ln y) G(y)dy = 1   &\text{ if } N=2.
\end{array}
\right.
\label{A15}
\end{equation}
It is then an easy exercice to extend $G$ on
$[\delta , \frac{1}{2}-\delta ]\cup [\frac{1}{2}+\delta, 1-\delta ]$
 in such a way  that $G$ is smooth and \eqref{A2}, \eqref{A12} and
\eqref{A15} are satisfied. This concludes the proof of Lemma \ref{lem2.1}.\qed

Let us turn now to the definition of the functions
$f_i$ and $g_i$ for $1\le i \le 3$.
Let $(t,z)$ ranges over $(-1,1)\times (0,1)$, so that
$(t,r)$ ranges over the domain
$$
\mathcal{O} :=\{ (t,r);\ -1<t<1,\ 0<r<\lambda (t) \}.
$$
Differentiating in \eqref{B1}, we obtain
$$
v_t=\sum_{i=0}^3
\big( \dot{f}_i g_i - \frac{\dot\lambda}{\lambda}
f_i z g_i^{(1)}\big),
$$
$$
v_{rr}+\frac{N-1}{r} v_r
=\sum_{i=0}^3 \lambda ^{-2}f_i
(g_i^{(2)} +\frac{N-1}{z} g_i^{(1)}),
$$
where ${\dot f}_i:= df_i/dt$ and $g_i^{(j)} := d^j g_i / dz^j$.
Let us introduce the function $\mathcal{V}=\mathcal{V}(t,z)$ defined by
\begin{eqnarray}
\mathcal{V} &:=& \lambda ^2 [v_{rr}+\frac{N-1}{r}v_r-v_t]\nonumber \\
&=& \sum_{i=0}^3 [f_i( g_i^{(2)}  +\frac{N-1}{z}g_i^{(1)})
+ z \lambda \dot\lambda
f_ig_i^{(1)} -\lambda ^2 \dot{f}_i g_i].
\label{A20}
\end{eqnarray}
We aim to define $f_i$ and $g_i$ so that
\begin{eqnarray*}
\mathcal{V}=\mathcal{V}_z=\mathcal{V}_{zz}=0 &\text{ for }& -1<t<1,\ z=\frac{1}{2},\\
\mathcal{V}_{zzz}\ge c f_0 &\text{ for }& -1<t<1,\
z=\frac{1}{2}
\end{eqnarray*}
for some constant $c>0$.
By \eqref{A1}, \eqref{A4} and (\ref{A20}),
\begin{eqnarray*}
\mathcal{V}(\cdot ,\frac{1}{2}) &=&
\frac{1}{2} \lambda \dot\lambda f_0g^{(1)}_0 (\frac{1}{2})
-\lambda ^2 \dot f_0 g_0(\frac{1}{2})\\
&&\quad
+ \sum_{i=1}^3 [f_i\big( g_i^{(2)} (\frac{1}{2})
+2(N-1)g_i^{(1)} (\frac{1}{2})\big)
+\frac{1}{2}\lambda \dot\lambda f_i g^{(1)}_i(\frac{1}{2})
-\lambda ^2 {\dot f}_i g_i(\frac{1}{2})].
\end{eqnarray*}
We impose the condition
\be
g_i^{(j)} (\frac{1}{2})=
\left\{
\begin{array}{ll}
1\quad &\text{if } i=1 \text{ and } j=2,\\
0 &\text{otherwise,}
\end{array}
\right.
  \text{ for } 1\le i\le 3,\ 0\le j\le 2.
\label{A21}
\ee
It follows that
$$
\mathcal{V}(\cdot ,\frac{1}{2}) =
\frac{1}{2}\lambda {\dot\lambda} f_0 g_0^{(1)}(\frac{1}{2})
-\lambda ^2 {\dot f}_0 g_0 (\frac{1}{2}) + f_1.
$$
The function $f_1$ is then defined by
\be
\label{A22}
f_1:= -\frac{1}{2}\lambda \dot\lambda f_0g_0^{(1)}
(\frac{1}{2}) + \lambda ^2 {\dot f}_0 g_0(\frac{1}{2}),
\ee
so that
\be
\label{A23}
\mathcal{V}(\cdot ,\frac{1}{2})=0 \text{ on } (-1,1).
\ee
Differentiating with respect to $z$ in \eqref{A20} yields, using once more (\ref{A4}),
\begin{eqnarray}
\mathcal{V}_z &=& f_0 G^{(1)} + z \lambda \dot\lambda
f_0 g_0^{(2)} + (\lambda \dot\lambda f_0 - \lambda ^2
{\dot f}_0) g_0^{(1)} \nonumber \\
&&+\sum_{i=1}^3
[f_i(g_i^{(3)} +\frac{N-1}{z}g_i^{(2)}
-\frac{N-1}{z^2} g_i^{(1)}) +z\lambda \dot\lambda
f_ig_i^{(2)} +(\lambda \dot\lambda f_i
-\lambda ^2 {\dot f}_i) g_i^{(1)}].
\label{A24}
\end{eqnarray}
We infer from \eqref{A1}, \eqref{A21} and \eq{A24} that
\begin{eqnarray*}
\mathcal{V}_z(\cdot ,\frac{1}{2})
=&\frac{1}{2} \lambda \dot \lambda f_0 g_0^{(2)}
(\frac{1}{2})
+(\lambda \dot\lambda f_0 -\lambda ^2 {\dot f}_0)
g_0^{(1)} (\frac{1}{2})
+\big( \sum_{i=1}^3 f_ig_i^{(3)}(\frac{1}{2})\big) +f_1(2(N-1)
+\frac{1}{2}\lambda \dot\lambda).
\end{eqnarray*}
We impose the condition
\be
g_i^{(3)}(\frac{1}{2})
=\left\{
\begin{array}{ll}
1   & \text{if } i=2,\\
0 &\text{if } i\in \{ 1,3\},
\end{array}
\right.
\label{A25}
\ee
and define $f_2$ as
\be
f_2 :=
-[  f_1 (2(N-1)  +\frac{1}{2}\lambda \dot \lambda)
+\frac{1}{2} \lambda \dot\lambda
f_0 g_0 ^{(2)} (\frac{1}{2})
+ (\lambda \dot\lambda f_0 -\lambda ^2 {\dot f}_0)
g_0^{(1)} (\frac{1}{2}) ].
\label{A26}
\ee
It follows that
\be
\mathcal{V}_z(\cdot ,\frac{1}{2}) =0   \text{ on } (-1,1).
\label{A27}
\ee
Differentiating (\ref{A24}) with respect to $z$, we get
\begin{eqnarray}
\mathcal{V}_{zz}
&=& f_0 G^{(2)} + z\lambda \dot\lambda f_0 g_0^{(3)}
+(2\lambda \dot\lambda f_0 - \lambda ^2 {\dot f}_0)g_0^{(2)}
\nonumber\\
&&\quad +\sum_{i=1}^3 [ f_i(g_i^{(4)}
+\frac{N-1}{z}g_i^{(3)} -2 \frac{N-1}{z^2} g_i^{(2)}
+2\frac{N-1}{z^3} g_i ^{(1)}) + z\lambda \dot\lambda
f_i g_i ^{(3)}\nonumber\\
&&\qquad\qquad  + (2\lambda \dot\lambda f_i -\lambda ^2 \dot f_i)
g_i^{(2)}], \label{A28}
\end{eqnarray}
which, together with \eqref{A1}, \eqref{A21} and \eqref{A25}, leads to
\begin{multline}
\mathcal{V}_{zz}(\cdot ,\frac{1}{2}) =
\frac{1}{2} \lambda \dot \lambda f_0g_0^{(3)}
(\frac{1}{2}) + (2\lambda \dot\lambda  f_0 -\lambda ^2
{\dot f}_0) g_0^{(2)}(\frac{1}{2}) \\
+ \big( \sum_{i=1}^3 f_i g_i ^{(4)}(\frac{1}{2}) \big)
+ 2(N-1) f_2
-8(N-1) f_1 + \frac{1}{2} \lambda \dot\lambda f_2
+(2\lambda \dot \lambda f_1 - \lambda ^2 {\dot f}_1).
\end{multline}
We impose the condition
\be
g_i^{(4)} (\frac{1}{2}) =
\left\{
\begin{array}{ll}
1   \text{ if } i=3,\\[3mm]
0   \text{ if } i\in \{1,2\} .
\end{array}
\right.
\label{A29}
\ee
and define $f_3$ as
\begin{eqnarray}
f_3&:=&-[
(2(N-1)+\frac{1}{2}\lambda \dot\lambda )f_2
+(2\lambda \dot\lambda -8(N-1))f_1 -\lambda ^2 {\dot f}_1\nonumber\\
&&\qquad +\frac{1}{2} \lambda \dot\lambda f_0g_0^{(3)} (\frac{1}{2})
+(2\lambda \dot\lambda f_0 -\lambda ^2 {\dot f}_0)
g_0^{(2)} (\frac{1}{2}) ].
\label{A30}
\end{eqnarray}
This gives
\be
\label{A31}
\mathcal{V}_{zz}(\cdot ,\frac{1}{2}) = 0 \text{ on } (-1,1).
\ee
By \eqref{A1} and \eqref{A28}, we have, for
$(t,z)\in (-1,1)\times (\frac{1}{2}-\delta,
\frac{1}{2}+\delta )$,
$$
\mathcal{V}_{zzz}=6f_0 +\mathcal{R},
$$
where
\begin{eqnarray}
\mathcal{R}&:=& z\lambda \dot\lambda f_0g_0^{(4)}
+(3\lambda \dot\lambda f_0 -\lambda ^2 {\dot f}_0 ) g_0^{(3)}
\nonumber \\
&&\quad +\sum_{i=1}^3 [ f_i \big(
g_i ^{(5)} + \frac{N-1}{z} g_i^{(4)}
-3\frac{N-1}{z^2} g_i ^{(3)}
+6\frac{N-1}{z^3}g_i^{(2)} -6\frac{N-1}{z^4} g_i^{(1)} \big) \nonumber\\
&&\qquad \qquad +z \lambda \dot\lambda f_i g_i ^{(4)}
+(3\lambda {\dot\lambda} f_i -\lambda ^2 {\dot f}_i)
g_i^{(3)} ].   \label{A32}
\end{eqnarray}
Let $C$ denote various constants independent of $\ep$, $t$,
and $z$, which may vary from line to line. We claim that
\be
\label{A33}
|\mathcal{R}| \le C\ep ^2 f_0  \text{ for }
(t,z)\in (-1,1)\times (\frac{1}{2} -\delta ,
\frac{1}{2} +\delta ).
\ee
First, we have that
$$
|\mathcal{R}| \le C\big( |\lambda\dot\lambda f_0| +
\lambda ^2 |{\dot f}_0|
+\sum_{i=1}^3 (|f_i|+|\lambda\dot\lambda f_i|
+\lambda ^2 |{\dot f}_i| \big) .
$$
Since
\begin{equation}
\label{zz1}
\lambda {\dot \lambda}f_0=-4\ep ^2 t (1-t^2)^3 f_0,\quad
\lambda ^2 {\dot f}_0 =-2\ep ^2 t (1-t^2)^2 f_0,
\end{equation}
one may write for each $i\in \{ 1 , ..., 3\}$
$$f_i(t)=\ep ^2p_i(t,\varepsilon )f_0(t),$$
where $p_i\in \R [t,\varepsilon ]$.
Therefore, there exists some constant $C>0$ such that
$$
|\lambda\dot{\lambda}f_0|
+\lambda^2|\dot{f}_0| + \sum_{i=1}^3(|f_i|+\lambda^2|{\dot f}_i|)
\le C \ep ^2 f_0,
$$
and \eqref{A33} follows. We infer that for
$\varepsilon$ small enough
$$
\mathcal{V}_{zzz} \ge (6-C\ep ^2) f_0 \ge f_0
\qquad\text{ for } (t,z)\in [-1,1] \times (\frac{1}{2} -\delta ,
\frac{1}{2}+\delta ).
$$
In view of the definitions of $\mathcal{V}$ and of $f_i$ for $1\le i\le 3$,
we may write
\begin{equation}\label{defA}
\mathcal{V}(t,z)=f_0(t)\sum_{j=1}^p
P_j(t)k_j(z) =: f_0(t) A(t,z),
\end{equation}
where $p\ge 1$, $P_j\in \R [t]$, $k_j\in C^\infty ([0, +\infty ))$.
Since
$$A(\cdot ,\frac{1}{2})=A_z(\cdot ,\frac{1}{2})
=A_{zz}(\cdot ,\frac{1}{2})=0$$
while
$$
A_{zzz}(t,z)\ge 1\quad \text{ for } \quad
(t,z)\in [-1,1]\times (\frac{1}{2}- \delta , \frac{1}{2}+\delta ),
$$
we conclude that we can write
\be
\label{A40}
A(t,z)=(z-\frac{1}{2})^3 \varphi(t,z)
  \text{ for }
t\in [-1,1],\ z\in (\frac{1}{2}-\delta , \frac{1}{2} +\delta ),
\ee
where
$\varphi \in C^\infty ( [-1,1]_t\times
(\frac{1}{2}-\delta , \frac{1}{2} +\delta )_z)$
and
$$
\varphi (t,z) > 0
  \text{ for } t\in [-1,1],\  z\in (\frac{1}{2}-\delta ,
\frac{1}{2} + \delta ).
$$
\, From \eqref{B4}, (\ref{A4}), (\ref{A20}) and (\ref{defA}) we have that, for
$t\in [-1,1]$ and $|z-\frac{1}{2}| \ge \frac{\delta}{2}$,
$$
A(t,z) = G + z \lambda {\dot \lambda} g_0^{(1)}
-\lambda ^2 \frac{{\dot f}_0}{f_0} g_0,
$$
which, combined to \eqref{zz1} and  \eqref{A2}, yields
$$
|A(t,z)-G | \le C\varepsilon ^2|G (z)|.
$$
It follows that for $\varepsilon >0$ small enough,
\be
\label{A41}
|A(t,z)| > \frac{1}{2} |G(z)| > 0
  \text{ for }
t\in [-1,1],\ z\in [0,\frac{1}{2} -\frac{\delta}{2}]\cup
[\frac{1}{2}+\frac{\delta}{2}, 1).
\ee
Gathering \eqref{A40}-\eqref{A41} we obtain that
$$
\lambda ^{-2}\mathcal{V}=\lambda ^{-2} f_0(t) A(t,z)=B(t,z)^3,\,
t\in (-1,1),\ z\in [0,1),
$$
for some $B\in C^\infty ((-1,1)\times [0,1))$. Define now
$V(t,x)$ by
$$
V(t,x):=v(t,|x|)= \sum_{i=0}^3
f_i(t) g_i( \frac{|x|}{\lambda (t)} ).
$$
\, From \eqref{B4} and \eqref{A6a} we have that
$$
V(t,x)=f_0(t)(1- \frac{|x|^2}{\lambda ^2(t)})
\text{ for }  |x|<\delta |\lambda (t)|.
$$
Combined with \eqref{B5}, \eqref{B6} and \eqref{A6b}, this yields
\be
\label{A50}
V\in C^\infty (\R\times \R ^N).
\ee
On the other hand, it follows from \eqref{B4}, \eqref{A7},
\eqref{A20}, (\ref{defA}), \eqref{A40} and \eqref{A41} that
$$
\lambda ^{-2}\mathcal{V}
=\lambda ^{-2}f_0 (z-\frac{1}{2}) ^3 (1-z^2)^{-4}
e^{-\frac{1}{1-z^2}} \psi (t,z) \text{ for } (t,z)\in (-1,1)\times [0,1),
$$
for some function $\psi \in C^\infty ([-1,1]\times [0,1])$ with
$|\psi (t,z)|\ge \eta >0$ on $[-1,1]\times [0,1]$. We observe that
the function
\be
\label{A51}
(\lambda ^{-2}\mathcal{V})^{\frac{1}{3}}
=(\lambda ^{-2}f_0)^{\frac{1}{3}}
(\frac{r}{\lambda} -\frac{1}{2})(1-\frac{r^2}{\lambda ^2})^{-\frac{4}{3}}
e^{-\frac{1}{3}\, \frac{1}{1-\frac{r^2}{\lambda ^2}}}
\psi ^{\frac{1}{3}}(t,\frac{r}{\lambda})
\ee
when extended by 0 for $|t|\ge 1$ or $r\ge \lambda (t)$, is of class
$C^\infty$ on $\R _t\times [0,+\infty )_r$. Therefore
$$
(\Delta V -V_t)^{\frac{1}{3}} (t,x) =(\lambda ^{-2} \mathcal{V})^{\frac{1}{3}}
(t,\frac{|x|}{\lambda (t)})
$$
is of class $C^\infty$ on $\R \times \R ^N \setminus
([-1,1]\times \{ 0\})$, and on a neighborhood of
$(-1,1)\times \{ 0\}$ by \eqref{A50} and the fact that
$$
(\Delta V -V_t)(t,0) = -2N \ep ^{-2} (1-t^2)^{-4}f_0(t)
+2t(1-t^2)^{-2} f_0(t) <0
$$
for $\ep >0$ small enough. The smoothness of $(v_t-\Delta v)^{1/3}$ near $(\pm 1,0)$ follows from
\eqref{B6} and  \eqref{A51}.
The proof of Theorem \ref{thm1} is complete. \qed

\section{Local null controllability around the trajectory $((\overline u, \overline v),\overline h)$}
\setcounter{equation}{0}
\label{seclocalnulbarubarv}

We consider the trajectory $((\overline u,\overline v),\overline h)$ of the control system \eq{system-0}
constructed in section \ref{sec-construction-baruetc}. Let $((u,v),h): Q\rightarrow \R^2\times \R$, and
let $((\zeta_1,\zeta_2),\widetilde h) :=((u-\overline u,v-\overline v),h-\overline h)$. Then $((u,v),h)$ is a trajectory of \eq{system-0}
if and only if $((\zeta_1,\zeta_2),\widetilde h)$ fulfills
\begin{equation}\label{system2}
    \left\{\begin{array}{ll}
    \displaystyle \zeta_{1,t}-\Delta \zeta_1=G_{11}(\zeta_1,\zeta_2)\zeta_1+G_{12}(\zeta_1,\zeta_2)\zeta_2
    +\widetilde h1_{\omega}&\hbox{ in }(0,T)\times \Omega, \\ \noalign{\medskip} \zeta_{2,t}-\Delta
    \zeta_2= G_{21}(\zeta_1,\zeta_2)\zeta_1 +G_{22}(\zeta_1,\zeta_2)\zeta_2
    &\hbox{ in } (0,T)\times\Omega,
    \\ \noalign{\medskip}
    \zeta_1=0,\quad \zeta_2=0&\hbox{ on }(0,T)\times\partial\Omega,
    \end{array}\right.
\end{equation}
where
\begin{eqnarray*}
\renewcommand{\arraystretch}{3}
\begin{array}{rcl}
G_{11}(\zeta_1,\zeta_2)(t,x)&:=&
\displaystyle
\int_0^1\frac{\partial
g}{\partial u}(\lambda \zeta_1(t,x)+\overline u(t,x),\zeta_2(t,x)+\overline v(t,x))\,d\lambda
\\
&=&
\left\{
\renewcommand{\arraystretch}{1}
\begin{array}{ll}
\displaystyle
\frac{g(\zeta_1+\overline u,\zeta_2+\overline
v)-g(\overline u,\zeta_2+\overline v)}{\zeta_1}(t,x) &\text{ if } \zeta_1(t,x)\not =0,
\\
\partial_u g(\overline u (t,x), \zeta_2(t,x) +\overline v(t,x) )&\text{ if } \zeta_1(t,x) =0,
\end{array}
\right.
\end{array}
\end{eqnarray*}
\begin{eqnarray*}
\renewcommand{\arraystretch}{3}
\begin{array}{rcl}
G_{12}(\zeta_1,\zeta_2)(t,x) &:=&
\displaystyle
\int_0^1\frac{\partial
g}{\partial v}(\overline u(t,x),\lambda \zeta_2(t,x)
+\overline v (t,x))\,d\lambda
\\
&=&
\medskip
\left\{
\renewcommand{\arraystretch}{1}
\begin{array}{ll}
\displaystyle
\frac{g(\overline u,\zeta_2+\overline
v)-g(\overline u,\overline v)}{\zeta_2}(t,x)&\text{ if } \zeta_2(t,x) \not =0,
\\
\partial_v g(\overline u(t,x), \overline v(t,x))&\text{ if } \zeta_2(t,x) =0,
\end{array}\right.
\end{array}
\end{eqnarray*}
\begin{gather*}
G_{21}(\zeta_1,\zeta_2):=(3\overline u ^2 +3 \overline u \zeta_1 +\zeta_1^2),
\\
\text{and}\qquad\qquad G_{22}(\zeta_1,\zeta_2):=R.
\qquad\qquad\qquad
\end{gather*}
Note that
\begin{gather}
\label{G21bar}
G_{21}(0,0)=3\overline{u}^2.
\end{gather}
By \eq{barunot0} and \eq{G21bar},
there exist $t_1\in(0,T)$, $t_2\in(0,T)$, a nonempty  open subset $\omega_0$
of $\Omega$ and $\overline M>0$ such that
\begin{gather}
\overline \omega_0\subset \omega, \, t_1<t_2,
\\
\label{G21>}
G_{21}(0,0)(t,x)\geqslant \frac{2}{\overline M}, \forall (t,x)\in  (t_1,t_2)\times \omega_0.
\end{gather}
Increasing $\overline M>0$ if necessary, we may also assume that
\begin{gather}
\|G_{ij}(0,0)\|_{L^\infty(Q)}\leqslant \frac{\overline M}{2},\, \forall (i,j)\in \{1,2\}^2.
\label{Gijleq}
\end{gather}

It is therefore natural to study the null controllability of the following
linear systems
\begin{equation}\label{linear2}
   \left\{\begin{array}{ll}
    \displaystyle \zeta_{1,t}-\Delta \zeta_1=a_{11}(t,x)\zeta_1+a_{12}(t,x)\zeta_2+h1_{\omega}&\hbox{ in }(0,T)\times \Omega, \\ \noalign{\medskip} \zeta_{2,t}-\Delta
    \zeta_2=a_{21}(t,x)\zeta_1+a_{22}(t,x)\zeta_2
        &\hbox{ in } (0,T)\times\Omega,\\ \noalign{\medskip}
         \zeta_1=0,\quad \zeta_2=0&\hbox{ on }(0,T)\times\partial\Omega,
    \end{array}\right.
\end{equation}
under the following assumptions
\begin{gather}
\label{a21>}
a_{21}(t,x)\geqslant \frac{1}{\overline M}, \forall (t,x)\in  (t_1,t_2)\times \omega_0.
\\
\label{aijbounded}
\|a_{ij}\|_{L^\infty(Q)}\leqslant \overline M,\, \forall (i,j)\in \{1,2\}^2.
\end{gather}
We will do this study in sub-section \ref{subseccontrolabbility-linearized}.
In sub-section \ref{sec-fixed-point},
we deduce from this study the local null controllability around the trajectory $((\overline u,\overline v),\overline h)$, and therefore get Theorem \ref{th}.

\subsection{Null controllability of a family of linear control systems}
\label{subseccontrolabbility-linearized}
Let $\mathcal{E}$ be the set of
$(a_{11}, a_{12}, a_{21}, a_{22})\in L^\infty (Q)^4$ such that
\eq{a21>} and \eq{aijbounded} hold. The goal of this sub-section is to prove the following lemma.

\begin{lem}
\label{leminf}
There exists $C>0$ such that, for every $(a_{11}, a_{12}, a_{21}, a_{22})\in \mathcal{E}$ and for every
$(\alpha_1,\alpha_2)\in L^2(\Omega)^2$, there exists a control
$h\in L^{\infty}(Q)$ satisfying
\begin{equation}
\label{regularcontrol}
\|h\|_{L^\infty(Q)}
\leq C(\|\alpha_1\|_{L^2(\Omega)}+\|\alpha_2\|_{L^2(\Omega)})
\end{equation}
such that the solution to the Cauchy problem
\begin{equation}
\label{linear2-Cauchy}
   \left\{\begin{array}{ll}
    \displaystyle \zeta_{1,t}-\Delta \zeta_1=a_{11}(t,x)\zeta_1+a_{12}(t,x)\zeta_2+h1_{\omega}&\hbox{ in }(0,T)\times \Omega, \\ \noalign{\medskip} \zeta_{2,t}-\Delta
    \zeta_2=a_{21}(t,x)\zeta_1+a_{22}(t,x)\zeta_2
        &\hbox{ in } (0,T)\times\Omega,\\ \noalign{\medskip}
         \zeta_1=0,\quad \zeta_2=0&\hbox{ on }(0,T)\times\partial\Omega,
    \\\noalign{\medskip}
    \zeta_1(0,\cdot)=\alpha_1,\, \zeta_2(0,\cdot)=\alpha_2&\text{ in } \Omega,
    \end{array}\right.
\end{equation}
satisfies
\begin{gather}
\label{v1Tv2T=0}
\zeta_1(T,\cdot)=0 \text{ and } \zeta_2(T,\cdot)=0.
\end{gather}
\end{lem}

Note that the coefficients $a_{jk}$ are in $L^{\infty}(Q)$, so that we have existence and uniqueness of solutions of (\ref{linear2-Cauchy}) in $C^0([0,T];L^2(\Omega)^2)\cap L^2(0,T;H^1(\Omega)^2)$. In order to prove Lemma~\ref{leminf}, we take $h(t,\cdot)=0$ for every $t\in (0,t_1)$. Note that  there then exists $C>0$ such that, for for every $(a_{11}, a_{12}, a_{21}, a_{22})\in \mathcal{E}$ and for every
$(\alpha_1,\alpha_2)\in L^2(\Omega)^2$, the solution to the Cauchy problem \eq{linear2-Cauchy} satisfies
$$
\|(\alpha_1^*,\alpha_2^*)\|_{L^2(\Omega)^2}\leqslant C \|(\alpha_1,\alpha_2)\|_{L^2(\Omega)^2},
$$
with
$$
(\alpha_1^*,\alpha_2^*):=(\zeta_{1}(t_1,\cdot),\zeta_{2}(t_1,\cdot)).
$$
 Then, our goal will be to find $h:(t_1,t_2)\times \Omega \rightarrow \R$ satisfying
 \begin{equation}
\label{regularcontrol-new}
\|h\|_{L^\infty((t_1,t_2)\times\Omega)}
\leq C(\|  \zeta_{1}(t_1,\cdot)\|_{L^2(\Omega)}+ \|\zeta_{2}(t_1,\cdot)\|_{L^2(\Omega)})
\end{equation}
such that the solution $(\zeta_1,\zeta_2)$ of the Cauchy problem
\begin{equation}
\label{linear2-Cauchy-t1t2}
   \left\{\begin{array}{ll}
    \displaystyle \zeta_{1,t}-\Delta \zeta_1=a_{11}(t,x)\zeta_1+a_{12}(t,x)\zeta_2+h1_{\omega}&\hbox{ in }(t_1,t_2)\times \Omega, \\ \noalign{\medskip} \zeta_{2,t}-\Delta
    \zeta_2=a_{21}(t,x)\zeta_1+a_{22}(t,x)\zeta_2
        &\hbox{ in } (t_1,t_2)\times\Omega,\\ \noalign{\medskip}
         \zeta_1=0,\quad \zeta_2=0&\hbox{ on }(t_1,t_2)\times\partial\Omega,
    \\\noalign{\medskip}
    \zeta_1(t_1,\cdot )=\alpha_1^*,\, \zeta_2(t_1,\cdot )=\alpha_2^*&\text{ in } \Omega,
    \end{array}\right.
\end{equation}
satisfies
\begin{gather}
\label{v1t2v2t2=0}
\zeta_1(t_2,\cdot)=0 \text{ and } \zeta_2(t_2,\cdot)=0.
\end{gather}
 Finally, we take $h(t,\cdot)=0$ for every $t\in (t_2,T)$. Of course, this construction provides
a control $h\in L^{\infty}((0,T)\times\Omega)$ driving the solution of (\ref{linear2-Cauchy}) to $(0,0)$ at time $t=T$. In paragraph~\ref{L2} we prove the existence of $h:(t_1,t_2)\times \Omega \rightarrow \R$ satisfying the required property but with  a $L^2$-bound instead of \eq{regularcontrol-new}. In paragraph \ref{Linf} we deal with the condition \eq{regularcontrol-new}. For the sake of simplicity, in these two paragraphs,  we write
 $(0,T)$ instead of $(t_1,t_2)$ and $Q$ instead of $(t_1,t_2)\times \Omega$.

\subsubsection{Controls in $L^2$}~\label{L2}
The goal of this paragraph is to prove a null controllability
result for the linear control systems \eq{linear2-Cauchy} with $L^2$ controls.
\begin{lem}
\label{lem}
There exists $C>0$ such that, for every $(a_{11}, a_{12}, a_{21}, a_{22})\in \mathcal{E}$ and for every
$(\alpha_1,\alpha_2)\in L^2(\Omega)^2$, there exists a control $h\in L^{2}(Q)$ satisfying
$$
\|h\|_{L^2(Q)}\leq C(\|\alpha_1\|_{L^2(\Omega)}+\|\alpha_2\|_{L^2(\Omega)})
$$
such that the solution to the Cauchy problem \eq{linear2-Cauchy} satisfies \eq{v1Tv2T=0}.
\end{lem}
In order to prove Lemma~\ref{lem}, we consider the associated adjoint
system
\begin{equation}\label{adjoint}
   \left\{\begin{array}{ll}
    \displaystyle -\varphi_{1,t}-\Delta\varphi_1=a_{11}(t,x)\varphi_1+a_{21}(t,x)\varphi_2&\hbox{ in }(0,T)\times \Omega,
    \\ \noalign{\medskip} -\varphi_{2,t}-\Delta
    \varphi_2=a_{12}(t,x)\varphi_1+a_{22}(t,x)\varphi_2
        &\hbox{ in } (0,T)\times\Omega,\\ \noalign{\medskip}
         \varphi_1=0,\quad \varphi_2=0&\hbox{ on }(0,T)\times\partial\Omega,
         \\ \noalign{\medskip}
\varphi_1(T,\cdot)=\varphi_{1,T}(\cdot),\quad
\varphi_2(T,\cdot)=\varphi_{2,T}(\cdot)&\hbox{ in }\Omega,
    \end{array}\right.
\end{equation}
and set $\varphi:= (\varphi_1,\varphi_2)$. For this system, we intend to prove the following observability
inequality:
\begin{equation}\label{OI}
\int_{\Omega}|\varphi(0,x)|^2\,dx\leq C\int\!\!\!
\int_{(0,T)\times \omega_0} |\varphi_1|^2\,dx\,dt.
\end{equation}
\, From estimate (\ref{OI}), with $C>0$ independent of $(a_{11}, a_{12}, a_{21}, a_{22})\in \mathcal{E}$, it is classical to deduce Lemma~\ref{lem}.

Let us recall the following  Carleman inequality, proved in  \cite[Chapter 1]{FurIma}, for the heat equation
with Dirichlet boundary conditions.
\begin{lem}
Let $\eta(t)=t^{-1}(T-t)^{-1}$. Let $\omega_1$ be a nonempty open set included in $\Omega$.
There exist a constant $C>0$ and a function $\rho\in C^2(\overline \Omega; (0,+\infty))$ such that, for every  $z\in H^1(0,T;L^2(\Omega))\cap L^2(0,T;H^2(\Omega))$ and for every $s\geq C$,

\begin{equation}\label{CI}
\begin{array}{l}\displaystyle
\int\!\!\!\int_{(0,T)\times\Omega}e^{-s\rho(x)\eta(t)}((s\eta)^3|z|^2+s\eta|\nabla
z|^2+(s\eta)^{-1}(|\Delta z|^2+|z_t|^2))\,dx\,dt
\\ \noalign{\medskip}\displaystyle
\leq
C\left(\int\!\!\!\int_{(0,T)\times\Omega}e^{-s\rho(x)\eta(t)}|z_t+\Delta
z|^2dx\,dt+\int\!\!\!\int_{(0,T)\times\omega_1}e^{-s\rho(x)\eta(t)}(s\eta)^3|z|^2\,dx\,dt\right).
\end{array}
\end{equation}
\end{lem}

For $\omega_1$ we take a nonempty open set of $\R^N$ whose closure (in $\R^N$) is included in $\omega_0$. Unless otherwise specified, we denote by $C$ various positive constants varying from
line to line which may depend of $\Omega$, $\omega_0$, $\omega_1$, $T$, $\rho$,
$\overline{M}$ and of other variables which will be specified later on. However they are independent  of $(a_{11}, a_{12}, a_{21}, a_{22})\in \mathcal{E}$,   of $(\varphi_{1,T},\varphi_{2,T})$, and of other variables which will
 be specified later on.

We start by applying (\ref{CI}) to $\varphi_1$ and $\varphi_2$
(solution of (\ref{adjoint})):
$$
\begin{array}{l}\displaystyle
\int\!\!\!\int_{(0,T)\times\Omega}e^{-s\rho(x)\eta(t)}((s\eta)^3|\varphi|^2+s\eta|\nabla
\varphi|^2+(s\eta)^{-1}(|\Delta \varphi|^2+|\varphi_t|^2))\,dx\,dt
\\ \noalign{\medskip}\displaystyle
\leq
C\left(\int\!\!\!\int_{(0,T)\times\Omega}e^{-s\rho(x)\eta(t)}(|a_{11}(t,x)\varphi_1+a_{21}(t,x)\varphi_2|^2
+|a_{12}(t,x)\varphi_1+a_{22}(t,x)\varphi_2|^2)dx\,dt \right.
\\ \noalign{\medskip}\displaystyle\left.
+\int\!\!\!\int_{(0,T)\times\omega_1}e^{-s\rho(x)\eta(t)}(s\eta)^3|\varphi|^2\,dx\,dt\right),
\end{array}
$$
for every $s\geqslant C$. Using \eq{aijbounded}
 and taking $s$ large enough, we have
\begin{equation}\label{Carle}
\begin{array}{l}\displaystyle
\int\!\!\!\int_{(0,T)\times\Omega}e^{-s\rho(x)\eta(t)}((s\eta)^3|\varphi|^2+s\eta|\nabla
\varphi|^2+(s\eta)^{-1}(|\Delta \varphi|^2+|\varphi_t|^2))\,dx\,dt
\\ \noalign{\medskip}\displaystyle
\leq C
\int\!\!\!\int_{(0,T)\times\omega_1}e^{-s\rho(x)\eta(t)}(s\eta)^3|\varphi|^2\,dx\,dt,
\end{array}
\end{equation}
for every $s\geqslant C$. Finally, we estimate the local integral of $\varphi_2$. We multiply the first equation in (\ref{adjoint}) by
$\chi(x)e^{-s\rho(x)\eta(t)}(s\eta)^3\varphi_2$, where
\begin{gather}
\label{propertieszeta}
\chi\in C^2(\overline \Omega; [0,+\infty)), \text{ the support of }
\chi \text{ is included in }\omega_0 \text{ and }
\chi =1 \text{ in }\omega_1.
\end{gather}
Integrating in
$(0,T)\times\omega_0$, this gives:
\begin{equation}\label{identity}
\begin{array}{l}\displaystyle
\int\!\!\!\int_{(0,T)\times\omega_0}a_{21}(t,x)\chi(x)e^{-s\rho(x)\eta(t)}(s\eta)^3|\varphi_2|^2\,dx\,dt
\\ \noalign{\medskip}\displaystyle
=
\int\!\!\!\int_{(0,T)\times\omega_0}\chi(x)e^{-s\rho(x)\eta(t)}(s\eta)^3\varphi_2(-\varphi_{1,t}-\Delta\varphi_1-a_{11}(t,x)\varphi_1)\,dx\,dt.
\end{array}
\end{equation}
Thanks to (\ref{a21>}) and \eq{propertieszeta}, the integral
in the left hand side of \eq{identity} is bounded from below by
$$
C^{-1}\int\!\!\!\int_{(0,T)\times\omega_1}e^{-s\rho(x)\eta(t)}(s\eta)^3|\varphi_2|^2\,dx\,dt.
$$
Let us now estimate the integral in the right hand side of
(\ref{identity}). Let $\varepsilon \in (0,1)$. From now on the constant $C>0$ may depend on $\varepsilon \in (0,1)$. Using \eq{aijbounded} (for $(i,j)=(1,1)$), we have that
\begin{equation}\label{esti1}
\begin{array}{l}\displaystyle
\left|\int\!\!\!\int_{(0,T)\times\omega_0}\chi(x)e^{-s\rho(x)\eta(t)}(s\eta)^3\varphi_2a_{11}(t,x)\varphi_1\,dx\,dt\right|
\\ \noalign{\medskip}\displaystyle
\leq \epsilon
\int\!\!\!\int_{(0,T)\times\Omega}e^{-s\rho(x)\eta(t)}(s\eta)^3|\varphi_2|^2\,dx\,dt
+C\int\!\!\!\int_{(0,T)\times\omega_0}e^{-s\rho(x)\eta(t)}(s\eta)^3|\varphi_1|^2\,dx\,dt.
\end{array}
\end{equation}
Next, for the time derivative term,
we integrate by parts with respect to $t$. We get
\begin{equation}\label{esti2}
\begin{array}{l}\displaystyle
-\int\!\!\!\int_{(0,T)\times\omega_0}\chi(x)e^{-s\rho(x)\eta(t)}(s\eta)^3\varphi_2\,\varphi_{1,t}\,dx\,dt
\\ \noalign{\medskip}\displaystyle
=
\int\!\!\!\int_{(0,T)\times \omega_0}\chi(x)e^{-s\rho(x)\eta(t)}(s\eta)^3\varphi_{2,t}\,\varphi_1\,dx\,dt
+\int\!\!\!\int_{(0,T)\times\omega_0}\chi(x)(e^{-s\rho(x)\eta(t)}(s\eta)^3)_t\varphi_2\,\varphi_1\,dx\,dt.
\end{array}
\end{equation}
Using that $|(e^{-s\rho(x)\eta(t)}(s\eta)^3)_t|\leq
Cs^4e^{-s\rho(x)\eta(t)}\eta(t)^5$ and Cauchy-Schwarz's
inequality, we can estimate this term in the following way
\begin{equation}\label{esti3}
\begin{array}{l}\displaystyle
\left|\int\!\!\!\int_{(0,T)\times\omega_0}\chi(x)e^{-s\rho(x)\eta(t)}(s\eta)^3\varphi_2\,\varphi_{1,t}\,dx\,dt\right|
\\ \noalign{\medskip}\displaystyle
\leq \epsilon
\int\!\!\!\int_{(0,T)\times\Omega}e^{-s\rho(x)\eta(t)}(s\eta)^{-1}(|\varphi_{2,t}|^2+(s\eta)^4|\varphi_2|^2)\,dx\,dt
\\ \noalign{\medskip}\displaystyle
+C\int\!\!\!\int_{(0,T)\times\omega_0}e^{-s\rho(x)\eta(t)}(s\eta)^7|\varphi_1|^2\,dx\,dt,
\end{array}
\end{equation}
for $s\geq C$. Finally, for the integral term with $\Delta\varphi _1$ we integrate by parts twice with respect to $x$ to get
$$
\begin{array}{l}\displaystyle
-\int\!\!\!\int_{(0,T)\times\omega_0}\chi(x)e^{-s\rho(x)\eta(t)}(s\eta)^3\varphi_2\,\Delta\varphi_1\,dx\,dt
=
-\int\!\!\!\int_{(0,T)\times\omega_0}\Delta(\chi(x)e^{-s\rho(x)\eta(t)}(s\eta)^3\varphi_2)\,\varphi_1\,dx\,dt.
\end{array}
$$
Using that
$$
|\Delta(\chi(x)e^{-s\rho(x)\eta(t)}(s\eta)^3\varphi_2)|\leq C
e^{-s\rho(x)\eta(t)}(s\eta)^3(|\Delta\varphi_2|+s\eta|\nabla\varphi_2|+(s\eta)^2|\varphi_2|), \,
(t,x)\in (0,T)\times\omega_0,
$$
in the previous identity together with Cauchy-Schwarz's inequality,
we deduce that
\begin{equation}\label{esti4}
\begin{array}{l}\displaystyle
-\int\!\!\!\int_{(0,T)\times\omega_0}\chi(x)e^{-s\rho(x)\eta(t)}(s\eta)^3\varphi_2\,\Delta\varphi_1\,dx\,dt
\\ \noalign{\medskip}\displaystyle
\leq \epsilon
\int\!\!\!\int_{(0,T)\times\Omega}e^{-s\rho(x)\eta(t)}(s\eta)^{-1}(|\Delta\varphi_{2}|^2+(s\eta)^2|\nabla\varphi_2|^2+(s\eta)^4|\varphi_2|^2)\,dx\,dt
\\ \noalign{\medskip}\displaystyle
+C\int\!\!\!\int_{(0,T)\times\omega_0}e^{-s\rho(x)\eta(t)}(s\eta)^7|\varphi_1|^2\,dx\,dt.
\end{array}
\end{equation}
Combining inequalities (\ref{esti1}), (\ref{esti3}) and
(\ref{esti4}) with (\ref{Carle}) and (\ref{identity}), and taking $\varepsilon \in (0,1)$ small enough, we obtain
\begin{equation}\label{Carle2}
\begin{array}{l}\displaystyle
\int\!\!\!\int_{(0,T)\times\Omega}e^{-s\rho(x)\eta(t)}((s\eta)^3|\varphi|^2+s\eta|\nabla
\varphi|^2+(s\eta)^{-1}(|\Delta \varphi|^2+|\varphi_t|^2))\,dx\,dt
\\ \noalign{\medskip}\displaystyle
\leq C
\int\!\!\!\int_{(0,T)\times\omega_0}e^{-s\rho(x)\eta(t)}(s\eta)^7|\varphi_1|^2\,dx\,dt,
\end{array}
\end{equation}
for every $s\geq C$.
From this estimate and taking into account the
dissipation (in time) of the heat system (\ref{adjoint}), one gets
for $s$ large enough
\begin{equation}
\label{Carle2-new-1}
\|\varphi(0,\cdot)\|_{L^2(\Omega)^2}^2
\leq C
\int\!\!\!\int_{(0,T)\times\omega_0}e^{-s\rho(x)\eta(t)}(s\eta)^7|\varphi_1|^2\,dx\,dt, \,
\end{equation}
which gives (\ref{OI}). Note that the constant $C$ in
\eqref{Carle2-new-1} may depend on $s$ at this time.
The proof of Lemma~\ref{lem} is finished.
\qed

\subsubsection{Controls in $L^{\infty}$}\label{Linf}

Let us remark that the proof in this sub-section follows ideas of
\cite{Barbu}. Let $\varepsilon \in (0,1)$. In this sub-section the constants $C>0$ do not depend on
$\varepsilon \in (0,1)$. We choose  $s>0$ large enough so that (\ref{Carle2})
(and therefore also \eq{Carle2-new-1}) holds.  Let $(\alpha_1,\alpha_2)\in L^2(\Omega)^2$. Let us consider, for each $\varepsilon>0$, the
extremal problem
\begin{equation}\label{minpb}
\inf_{h\in L^2((0,T)\times\omega_0)}\frac{1}{2} \int\!\!\!\int_{(0,T)\times\omega_0}e^{s\rho(x)\eta(t)}(s\eta)^{-7}|h|^2\,dx\,dt+\frac{1}{2\varepsilon}
\|\zeta(T,\cdot)\|^2_{L^2(\Omega)^2},
\end{equation}
where
$\zeta:=(\zeta_1,\zeta_2)$ is the solution of (\ref{linear2}) satisfying the initial condition
\begin{gather}
\label{initzeta}
\zeta_1(0,\cdot)=\alpha_1,\, \zeta_2(0,\cdot)=\alpha_2.
\end{gather}

We clearly have that there exists a (unique) solution of
(\ref{minpb}) $h^{\varepsilon}$ with
$(e^{s\rho(x)\eta(t)}(s\eta)^{-7})^{1/2}h^{\varepsilon}$ belonging
to $L^2((0,T)\times\omega_0)$. We extend $h^{\varepsilon}$ to all of $Q$ by letting
$h^{\varepsilon}:=0$ in $(0,T)\times (\Omega\setminus\omega_0)$. Let us call $\zeta^{\varepsilon}:=(\zeta_1^{\varepsilon},\zeta_2^{\varepsilon})$ the
solution of (\ref{linear2}) associated to $h^{\varepsilon}$ with, again, the initial condition \eq{initzeta}. The
necessary condition of minimum yields
\begin{equation}\label{necesscond}
\int\!\!\!\int_{(0,T)\times\omega_0}e^{s\rho(x)\eta(t)}(s\eta)^{-7}h^{\varepsilon}\, h\,dx\,dt
+\frac{1}{\varepsilon}\int_{\Omega}\zeta^{\varepsilon}(T)\cdot
\zeta(T)\,dx=0\quad \forall h\in L^2((0,T)\times\omega_0),
\end{equation}
where $\zeta:=(\zeta_1,\zeta_2)$ is the solution of
\begin{equation}\label{homocontrol}
\left\{\begin{array}{ll} \zeta_{1,t}-\Delta \zeta_1
=a_{11}(t,x)\zeta_1+a_{12}(t,x)\zeta_2+h1_{\omega_0}&\hbox{ in
}(0,T)\times\Omega,
\\ \noalign{\medskip}
\zeta_{2,t}-\Delta \zeta_2 =a_{21}(t,x)\zeta_1+a_{22}(t,x)\zeta_2&\hbox{ in
}(0,T)\times\Omega,
\\ \noalign{\medskip}
\zeta_1=0,\,\,\zeta_2=0&\hbox{ on } (0,T)\times\partial\Omega,
\\ \noalign{\medskip}
\zeta_1(0,\cdot)=\zeta_2(0,\cdot)=0&\hbox{ in }\Omega.
\end{array}
\right.
\end{equation}
Let us now introduce
$(\varphi^\varepsilon_1,\varphi^{\varepsilon}_2)$ the solution of
the following homogeneous adjoint system:
\begin{gather}
\label{eqphiepsilon}
\left\{\begin{array}{ll} -\varphi^{\varepsilon}_{1,t}
-\Delta\varphi^{\varepsilon}_1=a_{11}(t,x)\varphi_1^{\varepsilon}+a_{21}(t,x)\varphi^{\varepsilon}_2&\hbox{
in }(0,T)\times\Omega,
\\ \noalign{\medskip}
-\varphi^{\varepsilon}_{2,t}
-\Delta\varphi^{\varepsilon}_2=a_{12}(t,x)\varphi_1^{\varepsilon}+a_{22}(t,x)\varphi^{\varepsilon}_2&\hbox{
in }(0,T)\times\Omega,
\\ \noalign{\medskip}
\varphi^{\varepsilon}_1=\varphi^{\varepsilon}_2=0&\hbox{ on
}(0,T)\times\partial\Omega,
\\ \noalign{\medskip}\displaystyle
\varphi^{\varepsilon}(T,\cdot)=-\frac{1}{\varepsilon}\zeta ^{\varepsilon}(T,\cdot)&\hbox{
in } \Omega.
\end{array}\right.
\end{gather}
Then, the duality properties between $\varphi^{\varepsilon}$ and $\zeta$
provides
$$
-\frac{1}{\varepsilon}\int_{\Omega}\zeta ^{\varepsilon}(T)\cdot \zeta (T)\,dx=
\int\!\!\!\int_{ (0,T) \times \omega_0 }h\,\varphi^{\varepsilon}_1\,dx\,dt,
$$
which, combined with (\ref{necesscond}), yields
$$
\int\!\!\!\int_{(0,T)\times \omega_0 }h\,\varphi^{\varepsilon}_1\,dx\,dt
=\int\!\!\!\int_{(0,T)\times \omega_0}e^{s\rho(x)\eta(t)}(s\eta)^{-7}h^{\varepsilon}\,h\,dx\,dt
\quad\forall h\in L^2((0,T)\times\omega_0).
$$
Consequently, we can identify $h^{\varepsilon}$:
\begin{equation}
\label{hepsilon}
h^{\varepsilon}=e^{-s\rho(x)\eta(t)}(s\eta)^{7}\varphi^{\varepsilon}_11_{\omega_0}.
\end{equation}
\, From the systems fulfilled by $\zeta^{\varepsilon}$ and
$\varphi^{\varepsilon}$ we find, using (\ref{hepsilon}),
$$
-\frac{1}{\varepsilon}\|\zeta^{\varepsilon}(T,\cdot)\|^2_{L^2(\Omega)^2}=
\int\!\!\!\int_{(0,T)\times \omega_0}e^{-s\rho(x)\eta(t)}(s\eta)^7|\varphi^{\varepsilon}_1|^2
\,dx\,dt+\int_{\Omega}\varphi^{\varepsilon}(0,\cdot)\cdot \alpha \,dx,
$$
 with $\alpha := (\alpha_1,\alpha_2)$. Inequality (\ref{Carle2-new-1}) used for
$\varphi^{\varepsilon}$ tells us that
$$
\|\varphi^{\varepsilon}(0,\cdot)\|^2_{L^2(\Omega)^2}\leq
C\int\!\!\!\int_{(0,T)\times \omega_0}
e^{-s\rho(x)\eta(t)}(s\eta)^7|\varphi^{\varepsilon}_1|^2\,dx\,dt,
$$
so, using once more (\ref{hepsilon}),
\begin{equation}
\label{firstestimate}
\frac{1}{\varepsilon}
\|\zeta ^{\varepsilon}(T,\cdot)\|^2_{L^2(\Omega)^2}+
\|(e^{s\rho(x)\eta(t)}(s\eta)^{-7})^{1/2}h^{\varepsilon}\|^2_{L^2((0,T)\times\omega_0)}\leq
C \|\alpha \|^2_{L^2(\Omega)^2}.
\end{equation}
Consequently, we deduce the existence of a control $h$ such that
$(e^{s\rho(x)\eta(t)}(s\eta)^{-7})^{1/2}h\in L^2((0,T)\times\omega_0)$
(whose corresponding solution we denote by $\zeta$) such that
$\zeta (T,\cdot)= 0$  and
\begin{gather}
\label{ineqh}
\|(e^{s\rho(x)\eta(t)}(s\eta)^{-7})^{1/2}h\|_{L^2((0,T)\times\omega_0)}\leq
C \|\alpha \|_{L^2(\Omega)^2}.
\end{gather}

Let us finally bound the $L^{\infty}$-norm of the
control $h^{\varepsilon}$. For this, we develop now a boot-strap argument.

$\bullet$ Let
\begin{gather*}
\psi^{\varepsilon,0}:=
e^{-s\rho(x)\eta(t)/2}(s\eta)^{-1/2}\varphi^{\varepsilon},
\\
\psi^{\varepsilon,1}:=e^{-s\rho(x)\eta(t)/2}(s\eta)^{-5/2}\varphi^{\varepsilon}= \frac{1}{(s\eta)^2}\psi^{\varepsilon,0}.
\end{gather*}
Let $\mathcal{L}(\R^2;\R^2)$ be the vector space of linear maps from $\R^2$ into itself.  Using \eq{eqphiepsilon},
one easily checks that
$\psi^{\varepsilon,1}$ fulfills a  backward heat system with homogeneous
Dirichlet boundary condition and final null condition of the following form:
\begin{gather}
\label{systparabolic}
\left\{\begin{array}{ll}
-\psi^{\varepsilon,1}_t-\Delta \psi^{\varepsilon,1}= d_{1}& \text{ in } (0,T)\times \Omega,
\\ \noalign{\medskip}
\psi^{\varepsilon,1}=0 &\text{ on }(0,T)\times \partial \Omega,
\\ \noalign{\medskip}
\psi^{\varepsilon,1}(T,\cdot)=0& \text{ in } \Omega,
\end{array}
\right.
\end{gather}
with
\begin{gather}
\label{defd1}
d_{1}(t,x)=A_1(t,x)\psi^{\varepsilon,0}
+(s\eta)^{-1}\nabla \psi^{\varepsilon,0}\cdot \nabla \rho,
\end{gather}
where $A_1\in L^{\infty}(Q;\mathcal{L}(\R^2;\R^2))$  satisfy (see in particular \eq{aijbounded})
\begin{gather}
\label{estimateA1}
\|A_1\|_{L^{\infty}(Q;\mathcal{L}(\R^2;\R^2))}\leq C.
\end{gather}
In \eq{defd1} and in the following, we use the notation
$$
(\nabla \theta \cdot \nabla \rho) (t,x):=
(\nabla \theta _1(t,x)\cdot \nabla \rho (x),\nabla \theta _2(t,x)\cdot \nabla \rho (x)), \, \text{ for } \theta=(\theta_1,\theta_2): Q\rightarrow \R^2.
$$
Thanks to (\ref{Carle2}), \eq{hepsilon}, \eq{firstestimate}, \eq{defd1} and \eq{estimateA1},
\begin{gather}
\label{psi0in}
d_1 \in L^2(Q)^2 \text{ and }
\|d_1 \|_{L^2(Q)^2} \leq C \|\alpha \|_{L^2(\Omega)^2}.
\end{gather}
 For $r\in [1,+\infty)$, let $X_{r}:=L^{r}(0,T;W^{2,r}(\Omega)^2)\cap
W^{1,r}(0,T;L^{r}(\Omega)^2)$. We denote by $\|\cdot\|_{X_{r}}$ its usual norm. Let $X_{\infty}:=L^{\infty}(0,T;W^{1,\infty}(\Omega))^2$. We denote by $\|\cdot\|_{X_{\infty}}$ the usual $L^\infty$-norm.
Let
\begin{gather}
\label{defp1}
p_1:=2.
\end{gather}
\, From \eq{systparabolic}, \eq{psi0in}, \eq {defp1}
and a standard parabolic regularity theorem,
 we have
\begin{gather}
\label{psi1in}
\psi^{\varepsilon,1}\in
X_{p_1},\,
\|\psi^{\varepsilon,1}\|_{X_{p_1}}\leq
C\|\alpha \|_{L^2(\Omega)^2}.
\end{gather}

$\bullet$ For $k\in \mathbb{N}\setminus\{0\}$, let
$$\psi^{\varepsilon,k}:=
e^{-s\rho(x)\eta(t)/2}(s\eta)^{-1/2-2k}\varphi^{\varepsilon}=
\frac{1}{(s\eta)^{2k}}\psi^{\varepsilon,0}.
$$
Let us define, by induction on $k$, a sequence
$(p_k)_{k\in \mathbb{N}\setminus\{0\}}$ of elements of $[2,+\infty]$ by
\begin{eqnarray*}
p_k:=
\left\{
\begin{array}{ll}
\displaystyle
\frac{(N+2)p_{k-1}}{N+2-p_{k-1}}&\text{ if } p_{k-1}< N+2,
\\
2p_{k-1} &\text{ if } p_{k-1}= N+2,
\\
+\infty & \text{ if } p_{k-1}> N+2.
\end{array}
\right.
\end{eqnarray*}
One easily checks that
\begin{itemize}
\item [-] If $N=2l$, with $l\in \mathbb{N}\setminus \{0\}$, one has
\begin{gather*}
p_k=2\frac{N+2}{N-2(k-2)} \text{ if } k<l+2,
\\
p_{l+2}= 2(N+2),
\\
p_k=+\infty, \, \forall k\geqslant l+3.
\end{gather*}
\item [-] If $N=2l+1$, with $l\in \mathbb{N}$, one has
\begin{gather*}
p_k=2\frac{N+2}{N-2(k-2)} \text{ if } k\leqslant l+2,
\\
p_k=+\infty, \, \forall k > l+2.
\end{gather*}
\end{itemize}
In particular
\begin{gather}
\label{pkinfty}
p_k=+\infty, \forall k \geq \frac{N}{2}+3.
\end{gather}
We now use an induction argument on $k$. We assume that
\begin{gather}
\label{regpsik-1}
\psi^{\varepsilon,k-1}\in X_{p_{k-1}}\text{ and }
\|\psi^{\varepsilon,k-1}\|_{X_{p_{k-1}}}\leqslant C \|\alpha \|_{L^2(\Omega)^2},
\end{gather}
(now $C$ is allowed to depend on $k$) and that $\psi^{\varepsilon,k-1}$
fulfills a  heat system  of the following form:
\begin{gather}
\label{systparabolic-k-1}
\left\{
\begin{array}{ll}
-\psi^{\varepsilon,k-1}_t-\Delta \psi^{\varepsilon,k-1}= d_{k-1} & \text{ in } (0,T)\times \Omega,\\ \noalign{\medskip}
\psi^{\varepsilon,k-1}=0 &\text{ on }(0,T)\times \partial \Omega,\\ \noalign{\medskip}
\psi^{\varepsilon,k-1}(T,\cdot)=0& \text{ in } \Omega,
\end{array}
\right.
\end{gather}
with
\begin{gather}
\label{defdk-1}
d_{k-1}(t,x)=A_{k-1}(t,x)\psi^{\varepsilon,k-2}
+(s\eta)^{-1}\nabla \psi^{\varepsilon,k-2}\cdot \nabla \rho,
\end{gather}
where $A_{k-1}\in L^{\infty}(Q;\mathcal{L}(\R^2;\R^2))$  satisfies
\begin{gather}
\label{estimateAk-1}
\|A_{k-1}\|_{L^{\infty}(Q;\mathcal{L}(\R^2;\R^2))}\leq C.
\end{gather}
Note that we have just proved  above that this induction
assumption holds for $k=2$. Using \eq{systparabolic-k-1} and \eq{defdk-1}, one gets that $\psi^{\varepsilon,k}$
fulfills the following backward heat system with homogeneous
Dirichlet boundary condition and final null condition:
\begin{gather}
\label{systparabolic-k}
\left\{\begin{array}{ll}
-\psi^{\varepsilon,k}_t-\Delta \psi^{\varepsilon,k}= d_{k}& \text{ in } (0,T)\times \Omega,
\\ \noalign{\medskip}
\psi^{\varepsilon,k}=0 &\text{ on }(0,T)\times \partial \Omega,
\\ \noalign{\medskip}
\psi^{\varepsilon,k}(T,\cdot)=0& \text{ in } \Omega,
\end{array}
\right.
\end{gather}
with
\begin{gather}
\label{defdk}
d_{k}(t,x)=A_k(t,x)\psi^{\varepsilon,k-1}
+(s\eta)^{-1}\nabla \psi^{\varepsilon,k-1}\cdot \nabla \rho,
\end{gather}
where $A_k: Q \rightarrow \mathcal{L}(\R^2;\R^2)$ is defined by
\begin{gather}
\label{defAk}
A_k:=A_{k-1}+2\frac{\eta_t}{s^2\eta^3}\text{Id},
\end{gather}
$\text{Id}$  denoting the identity map of $\R^2$. From \eq{estimateAk-1} and \eq{defAk}, one gets that
\begin{gather}
\label{estimateAk}
A_k\in L^{\infty}(Q;\mathcal{L}(\R^2;\R^2)) \text{ and }\|A_k\|_{L^{\infty}(Q;\mathcal{L}(\R^2;\R^2))}\leq C.
\end{gather}
Let us recall the following embeddings between Sobolev spaces (see, e.g., \cite[Lemma 3.3, p. 80]{Lady}).
\begin{lem}
\label{lemSobolev}
Let $p\in (1,+\infty)$.
\begin{enumerate}
\item If $p<N+2$, let
$$
r:= \frac{(N+2)p}{N+2-p}.
$$
Then
$X_p$ is continuously embedded
in $L^r(0,T;W^{1,r}(\Omega)^2)$.
\item  If $p=N+2$, for every $r\in [1,+\infty)$,
$X_p$  is continuously embedded
in $L^r(0,T;W^{1,r}(\Omega)^2)$.
\item If $p>N+2$, $X_p$ is continuously embedded
in $L^\infty(0,T;W^{1,\infty}(\Omega)^2)$.
\end{enumerate}
\end{lem}

Applying Lemma \ref{lemSobolev} with $p=p_{k-1}$ and using \eq{regpsik-1},  we get that
\begin{gather}
\label{regpsikLk}
\psi^{\varepsilon,k-1}\in L^{p_{k}}(0,T;W^{1,p_{k}}(\Omega)^2)\text{ and }
\|\psi^{\varepsilon,k-1}\|_{L^{p_{k}}(0,T;W^{1,p_{k}}(\Omega)^2)}\leqslant C \|\alpha \|_{L^2(\Omega)^2}.
\end{gather}
\, From \eq{defdk}, \eq{estimateAk} and \eq{regpsikLk},
we have
\begin{gather}
\label{etimpsik}
d_{k}\in L^{p_{k}}(Q)^2 \text{ and }
\|d_{k}\|_{L^{p_{k}}(Q)^2}\leq C\|\alpha \|_{L^2(\Omega)^2}.
\end{gather}
 Using \eq{systparabolic-k}, \eq{etimpsik} and a classical parabolic regularity theorem
 (see, e.g., \cite[Theorem 9.1 p. 341--342]{Lady}, and (iii) of Lemma \ref{lemSobolev} if $p_k=+\infty$), we have
\begin{gather}
\label{regpsik}
\psi^{\varepsilon,k}\in X_{p_k}\text{ and }
\|\psi^{\varepsilon,k}\|_{X_{p_k}}\leq C\|\alpha \|_{L^2(\Omega)^2}.
\end{gather}
Hence \eq{regpsik} holds for every positive integer $k$. Let us choose an integer $k$ such that $k\geqslant (N/2)+3$. Then using \eq{pkinfty} and \eq{regpsik} we get that
\begin{gather}
\label{regpsikinfty}
\psi^{\varepsilon,k} \in L^\infty(Q)^2 \text{ and }
\|\psi^{\varepsilon,k}\|_{L^\infty (Q)^2}\leq C\|\alpha \|_{L^2(\Omega)^2}.
\end{gather}

This shows that the control defined in (\ref{hepsilon}) is bounded
in $L^{\infty}(Q)$ independently of $\varepsilon$ by $C\|\alpha \|_{L^2(\Omega)^2}$. This concludes the proof of
Lemma~\ref{leminf}.

\subsection{Local null controllability around the trajectory $((\overline u, \overline v), \overline h)$}
\label{sec-fixed-point}
Let $\nu >0$ be small enough so that, for every $z=(z_1,z_2)\in L^\infty(Q)^2$,
\begin{gather}
\label{propertynu}
\left(\|z\|_{L^\infty(Q)^2}\leq \nu\right)
\Rightarrow
\left((G_{11}(z_1,z_2),G_{12}(z_1,z_2),G_{21}(z_1,z_2),G_{22}(z_1,z_2))\in \mathcal{E}\right).
\end{gather}
(The existence of such a $\nu>0$ follows from \eq{G21>} and \eq{Gijleq}.) Let $\mathcal{Z}$ be the set of $z=(z_1,z_2)\in L^\infty(Q)^2$ such that $\|z\|_{L^\infty(Q)^2}\leq \nu$. By \eq{propertynu} and Lemma \ref{leminf}, there exists $C_0>0$ such that, for every $z=(z_1,z_2)\in \mathcal{Z}$ and for every $( \alpha _1 , \alpha _2)\in L^\infty(\Omega)^2$, there exists a
control $h\in L^{\infty}(Q)$ satisfying
\begin{equation}
\label{regularcontrolC0}
\|h\|_{L^\infty(Q)}
\leq C_0(\|\alpha_1\|_{L^2(\Omega)}+\|\alpha_2\|_{L^2(\Omega)}),
\end{equation}
such that the solution $(\zeta_1,\zeta_2)$ to the Cauchy problem
\begin{equation}
\label{linear2-Cauchy-new}
   \left\{\begin{array}{ll}
    \displaystyle \zeta_{1,t}-\Delta \zeta_1=G_{11}(z_1,z_2)\zeta_1+G_{12}(z_1,z_2)\zeta_2+h1_{\omega}&\hbox{ in }(0,T)\times \Omega, \\ \noalign{\medskip} \zeta_{2,t}-\Delta
    \zeta_2=G_{21}(z_1,z_2)\zeta_1+G_{22}(z_1,z_2)\zeta_2
        &\hbox{ in } (0,T)\times\Omega,\\ \noalign{\medskip}
         \zeta_1=0,\quad \zeta_2=0&\hbox{ on }(0,T)\times\partial\Omega,
    \\\noalign{\medskip}
    \zeta_1(0,\cdot)=\alpha _1,\, \zeta_2(0,\cdot)=\alpha _2&\text{ in } \Omega,
    \end{array}\right.
\end{equation}
satisfies \eq{v1Tv2T=0}.
We now define a set-valued mapping $B: \mathcal{Z} \rightarrow L^\infty(Q)^2$ as follows. Fix first any $(\alpha _1,\alpha _2)\in L^\infty (\Omega )^2$.
For any $z=(z_1,z_2)\in \mathcal{Z}$,  $B(z)$ is the set of $(\zeta_1,\zeta_2)\in L^\infty(Q)^2$ such that, for some $h\in L^\infty(Q)$ fulfilling
\eqref{regularcontrolC0},  $(\zeta_1,\zeta_2)$ is the solution of \eq{linear2-Cauchy-new} and this solution satisfies \eq{v1Tv2T=0}. As we have just pointed out, $B(z)$ is never empty.
Theorem \ref{th} will be proved if one can check that the set-valued mapping
$z\mapsto B(z)
$
has a fixed point (i.e. a point $z$ such that $z\in B(z)$) taking profit of the additional hypothesis
$\|(\alpha _1,\alpha _2)\|_{L^\infty(\Omega)^2}$  is small enough.
To get the existence of this fixed point, we apply
Kakutani's fixed point theorem (see, e.g., \cite[Theorem 9.B, page 452]{1986-Zeidler}): if
\begin{itemize}
\item [(i)] for every $z\in \mathcal{Z}$, $B(z)$ is a nonempty closed convex subset of $L^\infty (Q)^2$;
\item [(ii)] there exists a convex compact set $K\subset \mathcal{Z}$ such that
\begin{gather}\label{B(z)subset}
B(z)\subset
K,\, \forall z\in \mathcal{Z};
\end{gather}
\item [(iii)] $B$ is upper semi-continuous in $L^\infty (Q)^2$, i.e., for every closed subset $\mathcal{A}$ of
$\mathcal{Z}$, $B^{-1}(\mathcal{A}):=\{z\in \mathcal{Z};\, B(z)\cap \mathcal{A}\not =\emptyset\} $ is closed
 (see, e.g., \cite[Definition 9.3, page 450]{1986-Zeidler});
\end{itemize}
then there exists $z\in Z$ such that $z\in B(z)$.

Clearly (i) holds. Let us prove that (ii) holds. By standard estimates and using (\ref{regularcontrolC0}), there exists $C_1>0$ such that
\begin{gather}
\label{propC1}
\|\zeta \|_{L^{\infty}(Q)^2}\leq C_1 (\|
\alpha _1\|_{L^\infty(\Omega)}+\|\alpha _2\|_{L^\infty(\Omega)}),
 \, \forall z\in \mathcal{Z}, \,
\forall \zeta \in B(z).
\end{gather}
\, From now on we assume that $(\alpha _1,\alpha _2)\in  L^\infty(\Omega)^2$
satisfies
\begin{gather}
\label{choicer}
\|\alpha _1 \|_{L^\infty(\Omega)}+\|\alpha _2\|_{L^\infty(\Omega)}\leq \frac{\nu}{C_1}.
\end{gather}
\, From \eq{propC1} and \eq{choicer}, one has
\begin{gather}
\label{Bzsubset}
B(z)\subset \mathcal{Z}, \, \forall z \in \mathcal{Z}.
\end{gather}
Let $(\lambda _1,\lambda _2)\in L^\infty(Q)^2$ be the solution to the following Cauchy problem:
\begin{equation}
\label{defu*v*}
   \left\{\begin{array}{ll}
    \displaystyle \lambda _{1,t}-\Delta \lambda _1=0&\hbox{ in }(0,T)\times \Omega, \\ \noalign{\medskip}
\lambda _{2,t}-\Delta    \lambda _2=0
        &\hbox{ in } (0,T)\times\Omega,\\ \noalign{\medskip}
         \lambda _1=0,\quad \lambda _2=0&\hbox{ on }(0,T)\times\partial\Omega,
    \\\noalign{\medskip}
    \lambda _1(0,\cdot)=\alpha _1,\,
\lambda _2(0,\cdot)=\alpha _2 &\text{ in } \Omega.
    \end{array}\right.
\end{equation}
Let $\zeta_1^*:=\zeta_1-\lambda _1$ and $\zeta_2^*:=\zeta_2-\lambda _2$. Then $(\zeta_1^*, \zeta_2^*)$ is the solution to the following Cauchy problem
\begin{equation}
\label{Cauchyzeta*}
   \left\{\begin{array}{ll}
    \displaystyle \zeta_{1,t}^*-\Delta \zeta_1^*=D_1&\hbox{ in }(0,T)\times \Omega, \\ \noalign{\medskip} \zeta_{2,t}^*-\Delta
    \zeta_2^*=D_2
        &\hbox{ in } (0,T)\times\Omega,\\ \noalign{\medskip}
         \zeta_1^*=0,\quad \zeta_2^*=0&\hbox{ on }(0,T)\times\partial\Omega,
    \\\noalign{\medskip}
    \zeta_1^*(0,\cdot)=0,\, \zeta_2^*(0,\cdot)=0&\text{ in } \Omega,
    \end{array}\right.
\end{equation}
with
\begin{gather}
\label{D1}
D_1:=G_{11}(z_1,z_2)\zeta_1+G_{12}(z_1,z_2)\zeta_2+h1_{\omega},
\\
\label{D2}
D_2:= G_{21}(z_1,z_2)\zeta_1+G_{22}(z_1,z_2)\zeta_2.
\end{gather}
 Note that there exists $C_2>0$ such that
\begin{gather}
\label{estDinfty}
\|D_1\|_{L^\infty(Q)}+ \|D_2\|_{L^\infty(Q)}\leq C_2, \, \forall z\in \mathcal{Z}, \, \forall \zeta \in B(z).
\end{gather}
\, From \eq{Cauchyzeta*}, \eq{estDinfty} and
 a classical parabolic regularity theorem
(see, e.g., \cite[Lemma 3.3 p. 80 and Theorem 9.1 p. 341--342]{Lady}), $\zeta ^*:=(\zeta_1^*,\zeta_2^*)\in C^0(\overline Q)^2$ and there exists $C_3>0$ such that, for every
$z\in \mathcal{Z}$ and for every $\zeta \in B(z)$,
\begin{gather}
\label{regz*}
|\zeta^*(t,x)-\zeta^*(t',x')|\leq C_3 (|t-t'|^{1/2}+ |x-x'|) \,
\, \forall (t,x)\in \overline Q, \, \forall (t',x') \in \overline Q.
\end{gather}
Let $K^*$ be the set of $\zeta^*=(\zeta_1^*,\zeta_2^*)\in C^0(\overline Q)^2$ such that \eq{regz*} holds.
Then $(\lambda _1,\lambda _2)+K^*$ is a compact convex subset of $L^\infty(Q)^2$ and
\begin{gather}
\label{inclusion}
B(z)\subset (\lambda _1,\lambda _2)+K^*, \, \forall z\in \mathcal{Z}.
\end{gather}
 Then $K:=((\lambda _1,\lambda _2)+K^*)\cap \mathcal{Z}$
is a convex compact subset of $\mathcal{Z}$ such that \eq{B(z)subset} holds.

Let us finally prove the upper semi-continuity of $B$. Let $\mathcal{A}$ be a closed subset of $\mathcal{Z}$. Let
$(z^k)_{k\in \mathbb{N}}$ be a sequence of elements
in $\mathcal{Z}$, let $(\zeta^k)_{k\in \mathbb{N}}$  be a sequence of elements
in $L^\infty (Q)^2$, and  let $z\in \mathcal{Z}$ be such that
\begin{gather}
\label{convzk}
z^k\rightarrow z \text{ in } L^\infty(Q)\text{ as } k\rightarrow +\infty,
\\
\label{zetakincal}
\zeta^k \in \mathcal{A}, \, \forall k\in\mathbb{N},
\\
\label{zetakinBzk}
\zeta^k \in B(z^k), \, \forall k\in\mathbb{N}.
\end{gather}
By \eq{zetakinBzk}, for every $k\in \N$ there exists a control
$h^k\in L^\infty(Q)$ satisfying
\begin{gather}
\label{esthk}
\|h^k\|_{L^\infty(Q)}
\leq C_0(\|\alpha _1 \|_{L^2(\Omega)}+\|\alpha _2\|_{L^2(\Omega)}),
\end{gather}
(see (\ref{regularcontrolC0})) such that  $\zeta^k=(\zeta_1^k,\zeta_2^k)$ is the
solution of the Cauchy problem
\begin{equation}
\label{linear2-Cauchy-zetak}
   \left\{\begin{array}{ll}
    \displaystyle \zeta^k_{1,t}-\Delta \zeta^k_1=G_{11}(z^k_1,z^k_2)\zeta^k_1
    +G_{12}(z^k_1,z^k_2)\zeta_2^k +h^k1_{\omega}&\hbox{ in }(0,T)\times \Omega, \\ \noalign{\medskip} \zeta_{2,t}^k-\Delta
    \zeta^k_2=G_{21}(z^k_1,z^k_2)\zeta_1^k+G_{22}(z^k_1,z^k_2)\zeta^k_2
        &\hbox{ in } (0,T)\times\Omega,\\ \noalign{\medskip}
         \zeta^k_1=0,\quad \zeta^k_2=0&\hbox{ on }(0,T)\times\partial\Omega,
    \\\noalign{\medskip}
    \zeta^k_1(0,\cdot)=\alpha _1,\, \zeta^k_2(0,\cdot)=\alpha _2&\text{ in }
\Omega ,
    \end{array}\right.
\end{equation}
and this solution satisfies
\begin{gather}
\label{zetakfinal}
\zeta^k_1(T,\cdot)=0,\, \zeta^k_2(T,\cdot)=0.
\end{gather}
\, From (ii) and \eq{esthk},  there exists a strictly increasing sequence
  $(k_l)_{l\in \mathbb{N}}$ of integers, $h\in L^\infty(Q)$ and $\zeta \in \mathcal{Z}$ such that
\begin{gather}
\label{hkllim}
h^{k_l}\rightharpoonup h \text{ for the weak-* topology on } L^\infty(Q) \text{ as } l\rightarrow +\infty,
\\
\label{zetakllim}
\zeta^{k_l} \rightarrow \zeta \text{ in } L^\infty(Q)^2\text{ as } l\rightarrow +\infty.
\end{gather}
Note that, since $\mathcal{A}$ is closed, \eq{zetakincal} and \eq{zetakllim} imply that $\zeta \in \mathcal{A}$. Hence, in order to prove (iii),
it suffices to check that
\begin{gather}
\label{zetain B(z)}
\zeta\in B(z).
\end{gather}
Letting $l\rightarrow +\infty$ in \eq{linear2-Cauchy-zetak} and \eq{zetakfinal}, and using \eq{convzk}, \eq{hkllim} and \eq{zetakllim}, we get \eq{v1Tv2T=0} and \eq{linear2-Cauchy-new}. (The two equalities in
\eq{v1Tv2T=0} and  the two last equalities of \eq{linear2-Cauchy-new} have to be understood as equalities in $L^2(\Omega)$, for $\zeta ^{k_l}\rightharpoonup
\zeta $ in $X_2$.)
Letting $l\rightarrow +\infty$ in \eq{esthk} and using \eq{hkllim}, we get \eq{regularcontrolC0}. Hence \eq{zetain B(z)} holds. This concludes the proof of (iii) and of Theorem \ref{th}. \qed

\section{Appendix: Sketch of the proof of Theorem \ref{th2}.}
As the proof is very similar to those of Theorem \ref{th}, we limit ourselves
to pointing out the only differences. First, Lemma \ref{lem2.1} should be
replaced by
\begin{lem}
\label{lem2.1bis}
There exists a function $G\in C^\infty ([0,+\infty );\C)$
such that
\begin{eqnarray}
G(z) = (z-\frac{1}{2})^2 &\text{ for }&
\frac{1}{2} -\delta <z<\frac{1}{2} +\delta, \label{A1bis}\\
\text{Im}\, G(z) <0  &\text{ for }& 0<z<\frac{1}{2} -\delta ,
 \label{A2bis}\\
\text{Re}\, G(z) > \text{Im} \, G(z) >0 &\text{ for }& \frac{1}{2} + \delta < z < 1 - \delta
\label{A3bis}
\end{eqnarray}
and such that the solution $g_0$ to the Cauchy problem
\begin{eqnarray}
&&{g_0}''(z) +\frac{N-1}{z}{g_0}'(z)= G(z),\quad z>0,
\label{A4bis}\\
&&g_0(1)=g_0'(1)=0, \label{A5bis}
\end{eqnarray}
satisfies
\begin{eqnarray}
g_0(z)&=& 1-z^2\text{ if }\  0<z<\delta, \label{A6abis}\\
g_0(z)&=& e^{-\frac{1}{1-z^2}}
\text{ if }\ 1-\delta < z <1, \label{A6bbis}\\
g_0(z)&=& 0 \text{ if }\ z\ge 1. \label{A6cbis}
\end{eqnarray}
\end{lem}
The proof is carried out in the same way as for Lemma \ref{lem2.1}. Note that the conditions \eqref{A12} and \eqref{A15} are easily satisfied thanks to the change
of sign of Re\, $G$ and Im\, $G$. Theorem \ref{thm1} is still true for some
functions $V:(t,x)\in\R\times \R ^N\mapsto V(t,x)\in \C$ and
$K:(t,x)\in\R\times \R ^N\mapsto K(t,x)\in \C$ when  \eqref{E2} is replaced by
$$
V_t=\Delta V+RV+V^2.
$$
In the proof, we consider the same functions $\lambda (t)$ and $f_0(t)$ and search
$v(t,r)$ in the form
$$
v(t,r)=\sum_{i=0}^2 f_i(t)g_i(z).
$$
$f_1,f_2,g_1,g_2$ are defined in the same way as for Theorem \ref{thm1},
so that
\begin{eqnarray*}
{\cal V}={\cal V}_x=0\quad &\text{for}& -1<t<1,\ z=\frac{1}{2}\\
{\cal V}_{xx}=2f_0+{\cal R} \quad &\text{for}& -1<t<1,\ |z-\frac{1}{2}|<\delta
\end{eqnarray*}
with $|{\cal R}|\le C\varepsilon ^2 f_0$.
Letting
$A(t,z)=f_0(t)^{-1}{\cal V}(t,z)$, we have that
$$A(t,z)=(z-\frac{1}{2})^2\varphi (t,z)\quad
\text{ for } t\in [-1,1],\ z\in (\frac{1}{2}-\delta , \frac{1}{2}+ \delta )$$
where $\varphi \in C^\infty ([-1,1]_t\times (\frac{1}{2} -\delta,
\frac{1}{2} +\delta))$ and
$$\text{Re} \ \varphi (t,x)>f_0(t)\quad\text{ for }
t\in [-1,1],\ z\in (\frac{1}{2} -\delta , \frac{1}{2}+\delta ). $$
On the other hand
$$
|A(t,z)-G| \le C\varepsilon ^2 |G(z)|
$$
for $-1\le t\le 1$ and $| z-\frac{1}{2} | \ge \delta/2$.
Using \eqref{A2bis}, \eqref{A3bis}, it is then clear that for $\varepsilon$ small
enough we have $A(t,z)\not\in i\R ^+$ for $-1\le t\le 1$,
$z\in (0,1)\setminus \{ \frac{1}{2}\}$.
Defining the square root as an analytic function on the complement
of $i\R ^+$, we see that
$\lambda ^{-2}{\cal V} =
\lambda ^{-2} f_0(t)A(t,z)=B(t,z)^2$
for some $B\in C^\infty ((-1,1)\times [0,1))$. The end of the construction of
$(V,K)$ is as in the proof of Theorem \ref{thm1}. In the study of the local null
controllability around the trajectory $((\bar u, \bar v), \bar h)$,
the functions $G_{ij}$ are defined in the same way, except
$$
G_{21}(\zeta _1,\zeta _2)=2\bar u +\zeta _1.
$$
Therefore, \eqref{G21>} and \eqref{a21>} have to be changed respectively into
\begin{eqnarray*}
|\text{Im}\ G_{21}(0,0)(t,x)| &\ge& \frac{2}{\overline{M}}
\qquad \forall (t,x)\in (t_1,t_2)\times \omega _0,\\
|\text{Im}\ a_{21}(t,x)|      &\ge& \frac{1}{\overline{M}}
\qquad \forall (t,x)\in (t_1,t_2)\times \omega _0,
\end{eqnarray*}
Note that the functions in the control systems are complex-valued, so
that we have to conjugate the coefficients in the right hand side of the adjoint system \eqref{adjoint}. To estimate the local integral of $\varphi_2$, we multiply
the first equation in \eqref{adjoint} by
$\chi (x) e^{-s\rho (x)\eta (t)} (s\eta )^3 \overline{\varphi _2}$
(where $\overline{\varphi _2}$ stands for the conjugate of $\varphi _2$),
and take the absolute value of the imaginary part of the left hand side
of \eqref{identity}. The remaining part of the proof is the same as for
Theorem 1.

\end{document}